\documentclass[12pt,leqno,fleqn,epsfig]{article}
\usepackage{amssymb, epsfig, amsmath, amsthm}
\usepackage{color}
\usepackage{mathrsfs}       

\renewcommand{\r}{ \textcolor{red} } 
\newcommand{\R}{\mathbb{R}}

\newcommand{\N}{\mathbb{N}}

\newcommand{\esssup}{\operatorname*{ess\,sup}}

\newcommand{\curl}{\operatorname{curl}}

\newcommand{\bb}{\begin{equation}}
\newcommand{\ee}{\end{equation}}
\newcommand{\bq}{\begin{eqnarray}}
\newcommand{\eq}{\end{eqnarray}}
\newcommand{\bqn}{\begin{eqnarray*}}
\newcommand{\eqn}{\end{eqnarray*}}

\newcommand{\var}{\varepsilon}

\newcommand{\intl}{\int\limits}

\newcommand{\Beweisende}{\rule{0.2cm}{0.2cm}}

\newcommand{\intmw}{{\int\hspace{-830000sp}-\!\!}}

\newcounter{secnum}

\newtheorem{thm}{Theorem}[section]
\newtheorem{cor}[thm]{Corollary}
\newtheorem{lem}[thm]{Lemma}

\theoremstyle{definition}

\newtheorem{rem}[thm]{Remark}

\title{ Localized non  blow-up criterion of the Beale-Kato-Majda type for the 3D  Euler equations} 
 
\author{Dongho Chae$^*$  and J\"{o}rg Wolf $^\dagger$\\
\ \\
 Department of Mathematics\\
Chung-Ang University\\
 Seoul 156-756, Republic of Korea\\
$*$e-mail: dchae@cau.ac.kr\\
$\dagger$e-mail: jwolf2603@cau.ac.kr}
\date{}
\begin{document}
\maketitle
\begin{abstract}
We prove a localized non blow-up theorem  of the Beale-Kato-Majda type  for the  solution of the 3D incompressible Euler equations.\\
\
  \noindent{\bf
keywords:} Euler equation, local blow-up criterion

\end{abstract}

\section{Introduction}
\label{sec:-1}
\setcounter{secnum}{\value{section} \setcounter{equation}{0}
\renewcommand{\theequation}{\mbox{\arabic{secnum}.\arabic{equation}}}}

We consider the  homogeneous incompressible Euler equations\cite{eul} describing the fluid flows in  $\R^3$  , 
\begin{equation}
\partial _t v  + (v \cdot \nabla ) v  =  - \nabla p,\qquad  \nabla \cdot  v =0 \qquad \text{in}\quad \Bbb R^3\times (0, +\infty),
\label{euler}
\end{equation}
where $v=(v_1 (x,t), v_2 (x,t), v_3(x,t))$ is the velocity of the fluid,  and $p=p(x,t)$ represents the pressure.
The local in time well-posedness in  the Sobolev space $W^{k,p} (\Bbb R^3)$, $k>3/p+1$, $1<p<+\infty$, for the Cauchy problem of the system \eqref{euler}  is well-known by the result of  Kato-Ponce\cite{kat}. The question of finite time singularity for such local in time classical solution
 is still an outstanding open problem(see e.g.\cite{maj, con, bar,cha, tao} for surveys  of  the problem and the related results).  We say a local in time smooth solution $v\in C([0, T); W^{k,p}(\Bbb R^3))$, $k>3/p+1$, $1<p<+\infty$,
does not blow up (or equivalently regular) at $t=T$ if 
\begin{equation}\label{noblow}
\limsup_{t\nearrow T} \|v(t)\|_{W^{k,p}(\Bbb R^3)}<+\infty.
\end{equation}
The well-known Beale-Kato-Majda criterion\cite{bea} shows 
that \eqref{noblow} is guaranteed if and only if 
\begin{equation}\label{bk}
\int_{0} ^{T} \|\omega (t)\|_{L^\infty(\Bbb R^3)} dt <+\infty, \qquad \omega =\nabla \times v.
\end{equation} 
(see also \cite{con1, den} for  the other criteria, which are of geometric type).  In \cite{koz},  in particular,   Kozono-Taniuchi improved \eqref{bk}, replacing $\|\omega(t)\|_{L^\infty(\Bbb R^3)}$ in \eqref{bk} by a weaker norm $\|\omega(t)\|_{BMO(\Bbb R^3)}$.
The aim of the present paper is to 
   localize both of the previous criteria in \cite{bea, koz}.  For this purpose we use 
local BMO space. For $r>0$ and $x\in \Bbb R^n$ we denote  $B(x,r)=\{ y\in \Bbb R^n\, |\, |x-y|<r\}$, and $B(r)= B(0,r)$ below. 
By $ BMO(B(r))$ we denote the space of all $ u\in L^1(B(r) )$ such that 
\[
| u|_{ BMO(B(r) )} =  \sup_{\substack{z \in B(r)\\  0< \rho \le 2r} }  \intmw_{B(z,\rho) \cap  B(r) } | u- u_{ B(z, \rho )\cap B(r) }| dy  < +\infty,
\]
where we used the following notation for the average of $u$ over $\Omega \subset \Bbb R^n$.
$$u_\Omega = \intmw_\Omega udx.
$$
The space $ BMO(B(r))$ will be equipped with the norm 
\[
\| u\|_{ BMO(B(r))} = | u|_{ BMO(B(r))} + r^{ -n} \| u\|_{ L^1(B(1))}.
\]
Note that $ BMO(B(r))$ is continuously  embedded into $ L^q(B(r))$ for all $ 1 \le q < +\infty$. Indeed, in view of \eqref{B.3} 
in the  appendix below it holds 
\[
\| u\|_{ L^p(B(r))} \le c r^{ \frac{n}{q}} \| u\|_{ BMO(B(r))}. 
\]

{\em For  simplicity in presentation  we  assume that the possible blow-up occurs at the space-time  origin $(0,0)$, 
and consider the system \eqref{euler} in the local space time domain $ B(\rho)\times (-\rho, 0)$  throughout the paper. }
 Our aim in this paper is  the poof of the following form of local regularity criterion.
\begin{thm}
\label{thm1.1}
Let $(v,p)\in C^1(B(\rho)\times(-\rho,0))$ be a solution to \eqref{euler} such that  $v\in   C([- \rho ,0); W^{2,\, q}(B(\rho )))\cap L^\infty (- \rho, 0; L^2(B(\rho )) ) $
for some $ q\in(3, +\infty)$. If $ v$ satisfies 

\begin{equation}
 \intl_{- \rho }^{0} | \omega (s)|_{ BMO(B(\rho ) )}  d s < +\infty, 
\label{bkm}
\end{equation}
then  there exists no blow-up in $B(\rho)\times \{ 0\}$, namely 
$$ \lim\sup_{t\to 0^-}\|v(t)\|_{W^{2,\, q}(B(r))} <+\infty\qquad \forall  r\in (0, \rho).
$$

\end{thm}  

\begin{rem}
\label{rem1.1} From the obvious inequality $|\omega(s)|_{BMO(B(\rho)}\leq 2 \| \omega (s)\|_{ L^\infty(B(\rho ) )} $ we see that one can replace \eqref{bkm}
by
$
 \intl_{- \rho }^{0} | \omega (s)|_{ L^\infty(B(\rho ) )}  d s < +\infty 
$
in the above theorem. Thus, it generalizes both the original Beale-Kato-Majda criterion\cite{bea} and its improved version 
by Kozono-Taniuchi\cite{koz}. Moreover,  Theorem \ref{thm1.1} also provides  substantial advantage over the  global criterions of \cite{bea, koz} in the computational
test of the blow-up(see e.g. \cite{ker, hou} and references therein) at a specific point in a domain, since  we only need to compute the vorticities at  points in 
a small neigborhood of that point, not at whole points in the region.
\end{rem} 
 \ \\
 The contents of the paper is the following. In Section 2 we prove a localized version of the logarithmic Sobolev inequality. This is done by introducing suitable extension 
 operator of functions defined on a ball to the whole domain of $\Bbb R^n$.  In Section 3 we establish  several multiplicative inequalities to be used later. These correspond  to  localized versions of the Calder\'on-Zygmund type inequality in $\Bbb R^n$, which enable us to estimate the gradient of velocity in terms of the vorticity  together with lower order integral of  velocity. In Section 4 we prove a local $L_t ^\infty L_x ^p$ estimate for the vorticity. In order to do this we first prove a localized version of the Kozono-Taniuchi inequality(see \cite{koz1} for the original global version). The vorticity estimate deduced in this section, combined with our assumption of local energy bound,  implies $v\in L_t ^\infty L_x^\infty$ locally, which is an important step for our proof of the main theorem.
 In Section 5, using the results of previous sections, we complete the  proof of Theorem \ref{thm1.1}. This last part of the proof  is based on 
 the two new ingredients. One is the transform of the Euler system into new equations, which is  similar to the one in our previous paper\cite{cha1}. 
 The other one is  use a new iteration scheme of the Gronwall type. The corresponding iteration lemma is proved in Appendix A.
 
 \section{Local version of logarithmic Sobolev's inequality}
 \label{sec:-8}
 \setcounter{secnum}{\value{section} \setcounter{equation}{0}
 \renewcommand{\theequation}{\mbox{\arabic{secnum}.\arabic{equation}}}}
Our aim in this section is  to prove the following local version of the logarithmic Sobolev 
inequality.
\begin{lem}
\label{lem8.2} Let $B(r)$ be a ball in $\Bbb R^n$ with the radius $r>0$.
For every $ u\in W^{1,\, q}(B(r))$, $n<  q < +\infty$,  the following inequality holds true 
\begin{equation}
\| u \|_{ L^\infty(B(r))} \le c(1+ \| u\|_{ BMO(B(r))}) \log \Big(e+ c\|\nabla  u\|_{ L^{ q}(B(r))}+ c  r^{ -1+ \frac{n}{q}- \frac{n}{2}} \| u\|_{ L^2(B(r))}\Big)
\label{7.12}
\end{equation}
with a constant $ c>0$ depending only on $ n$ and $q$. 
\end{lem}

In order to prove the above lemma  we construct an extension operator, which is bounded with respect to both the $ BMO$ norm 
 and the Sobolev norm.
 \vspace{0.5cm}  
\begin{lem}
\label{lem8.1}
Let $ u \in W^{1,\, q}(B(r))\cap BMO(B(r))$ with $B(r)\subset \Bbb R^n$.  Let $ \phi \in C^{\infty}_{\rm c}(B(3r))$ denote a cut function such that $ 0 \le \phi \le 1$ 
in $ \R^{n}$, $ \phi \equiv 1$ on $ B(2r)$ and $ | \nabla \phi |^2 + | \nabla ^2 \phi | \le c r^{ -2}$. We define an extension $U$ of $u$ as follows.
\begin{equation}
U (x):=\begin{cases}
u(x) \quad  \text{ if}\quad x\in B(r) 
\\[0.3cm]
u (T(x))\phi (x)\quad \text{ if}\quad x\in B(r)^c,
\end{cases}
\label{8.1}
\end{equation}
where $ T: \R^{n}  \setminus \{ 0\} \rightarrow \R^{n}  \setminus \{ 0\}$ stands for the inversion with respect to $\partial B(r)$
\[
T(y) = \frac{r^2}{| y|^2} y,\quad  y\in \R^{n}  \setminus \{ 0\}. 
\] 

Then, $ U\in W^{1,\, q}(\R^{n})\cap BMO $. In addition, the following estimates hold true
\begin{align}
\| \nabla U \|_{L^{ q}} &\le c  (  \| \nabla  u \|_{ L^{ q}(B(r))}+ r^{ -1+ \frac{n}{q}- \frac{n}{2}}\|   u \|_{ L^{ 2}(B(r))}),
\label{8.2}
\\
| U| _{ BMO} &\le  c\| u\|_{ BMO(B(r))}. 
\label{8.2a}
\end{align}

\end{lem}

We first prove Lemma \ref{lem8.2}  assuming Lemma \ref{lem8.1} is true.\\

{\bf Proof of Lemma \ref{lem8.2}}: We denote by $ U $ the extension introduced in Lemma\,\ref{lem8.1}.  According to Lemma\,\ref{lem8.1} 
we get $ U \in W^{1,\, q}(\R^{n})$. In view of the logarithmic Sobolev embedding we infer 
\begin{equation}
\| u\|_{ L^\infty(B(r))} \le \| U \|_{\infty} \le 
c(1+ \| U \|_{ BMO} ) \log (e+ \|\nabla  U \|_{q}).
\label{7.13}
\end{equation}
Estimating the right-hand side of \eqref{7.13} by means of \eqref{8.2} and \eqref{8.2a}, the assertion follows.  \hfill \Beweisende    \\

{\bf Proof of Lemma \ref{lem8.1}}: 
First let us provide some basic properties of the map $ T$.  We compute 
\[
\partial _j T_i(y) = \frac{r^2}{| y|^2} \delta _{ ij} - 2r^2\frac{ y_iy_j}{| y|^4},\quad  y\in \R^{n}  \setminus \{ 0\},\quad i,j=1, \ldots, n. 
\]
Furthermore we get for $ k,l\in \{ 1, \ldots, n\}$
\begin{align*}
\partial _k T (y)\cdot \partial _l T (y)&= \frac{r^4}{| y|^4} (\delta _{ik} - 2\frac{ y_iy_k}{| y|^2}) (\delta _{il} - 2\frac{ y_iy_l}{| y|^2})
\\
&= \frac{r^4}{| y|^4} (\delta _{ kl} - 4\frac{ y_ly_k}{| y|^2}+ 4\frac{ y_ly_k}{| y|^2} )= \frac{r^4}{| y|^4} \delta _{ kl}. 
\end{align*}
Accordingly,
\[
(\det DT(y))^2 =\det (DT(y)^{ \top} DT(y)) = \det \Big(\frac{r^4}{| y|^4}I\Big) = \Big(\frac{r}{| y|}\Big)^{ 4n}.   
\]
This show that 
\begin{equation}
| \det DT(y)| = \Big(\frac{r}{| y|}\Big)^{ 2n}.   
\label{8.3}
\end{equation}
Next, we show that for every $ x,y\in \R^{n}  \setminus \{ 0\}$ it holds
\begin{equation}
| T(x)- T(y)| = r^2 \frac{| x-y|}{| x| | y|}. 
\label{8.4}
\end{equation}
Indeed, by an elementary calculus we find 
\begin{align*}
| T(x)- T(y)|^2  &=  \Big|\frac{r^2}{| x|^2} x-  \frac{r^2}{| y|^2} y\Big|^2
\\
&= \frac{r^4}{| x|^2 | y|^2} \Big|\frac{ | y|}{| x|} x-  \frac{| x|}{| y|} y\Big|^2 = 
\frac{r^4}{| x|^2 | y|^2} | x-y|^2.
\end{align*}
Whence, \eqref{8.4}. 
 
 \hspace{0.5cm}
Let $ B(x_0, \rho ) \subset \R^{n}$, $ 0< \rho <  r$ be any ball.  We discuss the three cases 
(i) $ B(x_0, \rho ) \subset B(r)$, (ii) $ B(x_0, \rho ) \subset B(r)^c$, and (iii) $ B(x_0, \rho ) \cap  \partial B(r ) \neq \emptyset$ separately.

\vspace{0.2cm}
(i) {\it The case $ B(x_0, \rho ) \subset B(r)$}: As $ U = u$ in $ B(r)$ it follows that  
\[
\intmw_{B(x_0, \rho )} | U (y) - U_{ B(x_0, \rho )} | dy =  
\intmw_{B(x_0, \rho )} | u (y) - u_{ B(x_0, \rho )} | dy \le  | u|_{ BMO(B(r))}. 
\]
(ii) {\it The case $ B(x_0, \rho ) \subset B(r)^c$}: By means of triangle inequality and Jensen's inequality we estimate  
\begin{align}
\intmw_{B(x_0, \rho )} | U  - U_{ B(x_0, \rho )} | dy \le \intmw_{B(x_0, \rho )}
\intmw_{B(x_0, \rho )} | \phi (y) u(T(y))-\phi (y') u(T(y'))| dydy'.  
\label{8.5} 
\end{align}
Using change of coordinates $ z = T (y), z'= T(y')$ ($ y= T(z)$, $ y' = T(z')$), and observing \eqref{8.3}, we obtain with the help of the transformation 
formula of the Lebesgue integral  
\begin{align}
 &\intmw_{B(x_0, \rho )}
\intmw_{B(x_0, \rho )} | \phi (y) u(T(y))-\phi (y') u(T(y'))| dydy'
\cr
&  \le \frac{1}{c_n^2 \rho ^{ 2n}}  \int_{T(B(x_0, \rho) )}
\int_{T(B(x_0, \rho ))} | \phi (T(z)) u(z)-\phi (T(z')) u(z')| \Big(\frac{r}{| z|}\Big)^{ 2n}\Big(\frac{r}{| z'|}\Big)^{ 2n} dzdz'.  
\label{8.6}
\end{align}
In case $ x_0 \in B(4r)^c$ it is readily seen that $ B(x_0, \rho ) \subset B(3r)^c$, and therefore $ U \equiv 0$ on $ B(x_0, \rho )$. 
Thus, we may assume $ x_0 \in B(4r)$. By means of \eqref{8.4} we infer for every $ y\in B(x_0, \rho ), $ since $|x_0|>r+\rho$ and $ |y|>r$,
\[
| T(x_0) - T(y)|  =  r^2 \frac{| x_0-y |}{| x_0|| y|}  \le | x_0-y|< \rho.  
\] 
This, together with $T(B(x_0, \rho ))  \subset B( r)$, implies that 
\begin{equation}
 T(B(x_0, \rho )) \subset B(T(x_0), \rho ) =B(T(x_0), \rho )\cap B(r).   
\label{8.7}
\end{equation}
With the help  of triangle inequality we obtain for $ z= T(y)$ and $y\in B(x_0, \rho )$
\begin{equation}
| z|= | T(y)|= \frac{r^2}{| y|}  \ge \frac{r^2}{ | y- x_0|+ | x_0|} >  \frac{ r^2}{ \rho +4r}  >   \frac{r}{5}  .   
\label{8.8}
\end{equation}
Now, the right-hand side of \eqref{8.6}  can be estimated by virtue of \eqref{8.7} and \eqref{8.8}. This gives 
\begin{align}
 &\intmw_{B(x_0, \rho )}
\intmw_{B(x_0, \rho )} | \phi (y) u(T(y))-\phi (y') u(T(y'))| dydy'
\cr
&  \le \frac{5^{ 4n}}{c_n^2 \rho ^{ 2n}}  \int_{T(B(x_0, \rho))}
\int_{T(B(x_0, \rho))} | \phi (T(z)) u(z)-\phi (T(z')) u(z')| dzdz'   
\cr
&  \le 5^{ 4n}  \intmw_{B(T(x_0), \rho)}
\intmw_{B(T(x_0), \rho )} |  u(z)-u(z')| dzdz'   
\cr
&\qquad + \frac{5^{ 4n}}{c_n^2 \rho ^{ 2n}}    \int_{T(B(x_0, \rho))}
\int_{T(B(x_0, \rho))} | (\phi (T(z)) -\phi (T(z'))) u(z')| dzdz'   = I+ II. 
\label{8.9} 
\end{align}  
Clearly,  $ I \le 5^{ 4n} | u|_{ BMO(B(r))}$. To estimate the second integral we first note that for all $ z \in T(B(x_0, \rho ))$  
\[
| \phi (T(z))-\phi (T(z'))| \le 2\| \nabla \phi \|_{ \infty} \rho \le c r^{ -1} \rho.  
\]
 This together with \eqref{8.7} shows that 
 \[
II \le c r^{ -1} \rho  \intmw_{B(T(x_0), \rho)}
 |  u(z')| dz' \le c \| u\|_{ BMO(B(r))}+ c r^{ -1} \rho | u_{ B(T(x_0), \rho )}|.    
\] 
Let $ k\in \N$ such that $ r_k=  2^{ -k}r \ge \rho$.  By means of 
triangle inequality and Jensen's inequality along with 
\[
|B(T(x_0), r_{ k}) \cap B(r)| \ge c r_{k}^{n}
\]
for some positive constant $c$ depending on $n$ only, we get 
\begin{align*}
& \Big || u_{ B(T(x_0), r_{ k}) \cap B(r)}| - | u_{ B(T(x_0), r_{ k-1}) \cap B(r)}|\Big | 
\\
&\quad \le  \Big | u_{ B(T(x_0), r_{ k}) \cap B(r)} -  u_{ B(T(x_0), r_{ k-1}) \cap B(r)}\Big | 
\\
&\quad \le  \intmw_{B(T(x_0), r_{ k}) \cap B(r)} \intmw_{B(T(x_0), r_{ k-1}) \cap B(r)}  |u(x)-u(y)| dxdy
\\
&\quad \le c r_{k}^{-2n} \int_{B(T(x_0), r_{ k}) \cap B(r)} \int_{B(T(x_0), r_{ k}) \cap B(r)}  |u(x)-u(y)| dxdy
\\
&\quad \le c  \intmw_{B(T(x_0), r_{ k}) \cap B(r)} \int_{B(T(x_0), r_{ k}) \cap B(r)}  |u(x)- u_{ B(T(x_0), r_{ k}) \cap B(r)}| dx \le | u|_{ BMO(B(r))}.
\end{align*}

Using the inequality above, along with  triangle inequality we find 
 \begin{align*}
& r_{ k}| u_{ B(T(x_0), r_{ k}) \cap B(r)}| 
\\
&\le r_{ k} \Big || u_{ B(T(x_0), r_{ k}) \cap B(r)}| - | u_{ B(T(x_0), r_{ k-1}) \cap B(r)}|\Big |  +
 r_{ k}| u_{ B(T(x_0), r_{ k-1}) \cap B(r)}| 
\\
&\le  2^{ -k}c r | u|_{ BMO(B(r))} + r_{ k-1}| u_{ B(T(x_0), r_{ k-1}) \cap B(r)}| . 
\end{align*}
Iterating the above inequality from $ k=N$ to $ k=1$, where $N$ is chosen so that  $ 2^{ -N-1}r < \rho \le 2^{ -N} r$, we get 
\begin{align*}
 \rho   | u_{ B(T(x_0), \rho )}| &\le c r | u|_{ BMO(B(r))} +  cr | u_{ B(T(x_0), r )\cap B(r)}| 
\\
&\le 
 cr  ( | u|_{ BMO(B(r))} + r^{ -n} \| u\|_{ L^1(B(r))}).   
\end{align*}
Accordingly, $ II \le c \| u\|_{ BMO(B(r))} $. Inserting the estimates of $ I$ and $ II$ into the 
right-hand side of \eqref{8.9}, we obtain from  \eqref{8.5}
\[
\intmw_{B(x_0, \rho )} | U  - U_{ B(x_0, \rho )} | dy \le c \| u\|_{ BMO(B(r))}. 
\]

(iii) {\it The case $ B(x_0, \rho ) \cap  \partial B(r ) \neq \emptyset$}: We decompose the integral as follows.
\begin{align*}
&\intmw_{B(x_0, \rho) }\intmw_{B(x_0, \rho )} |U  (y)-U  (y')| dy dy' 
\\
&\quad = \frac{1}{c_n^2 \rho ^{ 2n} }\int_{B(x_0, \rho)\cap B(r) }\int_{B(x_0, \rho )\cap B(r)} |u (y)-u(y')| dydy' 
\\
&\qquad +\frac{2}{c_n^2 \rho ^{ 2n} }\int_{B(x_0, \rho)\cap B(r) }\int_{B(x_0, \rho )  \setminus  B(r)} |u(T(y))-u  (y')| dy dy' 
\\
&\qquad +\frac{1}{c_n^2 \rho ^{ 2n} }\int_{B(x_0, \rho)  \setminus  B(r) }\int_{B(x_0, \rho )  \setminus  B(r)} |u(T(y))- u(T(y'))| dy dy' =I+II+III.
\end{align*}
We set  $ \hat{x} := \frac{r}{| x_0|} x_0\in \partial B(r)$. Since $ B(x, \rho )\cap B(r) \neq \emptyset$ it  must hold 
$ | x_0| < r+\rho $,  and therefore $ \hat{x} \in B(x_0, \rho )$. In particular, $ | x_0- \hat{x}| < \rho $.   
This shows that 
\[
| x_0- \hat{x}|^2 = | x_0|^2 \r{-} 2r| x_0|+ r^2 = (| x_0|-r)^2 < \rho ^2,
\]
which in turn implies 
\begin{equation}
\Big|| x_0|-r\Big| < \rho. 
\label{8.10}
\end{equation}
 By the triangle inequality we infer for every $ y\in B( x_0 , \rho )$
\[
| y- \hat{x} | \le | y-x_0| + | x_0- \hat{x}| < 2\rho. 
\]
For  $ y \in B(x, \rho )  \setminus B(r)$, noting that $ \hat{x}=T(\hat{x})$, and from  \eqref{8.4}  together with the facts  $ |y|>r$,  $|\hat{x}|= r$ we find 
\[
| T(y)- \hat{x} | = | T(y)- T(\hat{x})| =r^2 \frac{| y-\hat{x}|}{| y|| \hat{x}|} < 2\rho. 
\]
Consequently,
\begin{equation}
B(x_0  , \rho )\cap B(r) \subset B(\hat{x}, 2\rho )\cap B(r),\quad  T(B(x_0  , 2\rho )  \setminus B(r)) \subset B(\hat{x}, 2\rho )\cap B(r).
\label{8.11}
\end{equation}
Using \eqref{8.11}, we obtain 
\[
I \le c \intmw_{B(\hat{x}, 2\rho)\cap B(r) }\intmw_{B(\hat{x}, 2\rho )\cap B(r)} |u (y)-u(y')| dydy' 
\le c \| u\|_{ BMO(B(r))}. 
\]
To estimate $ II$ we proceed as in (ii).  For  $ y\in B(x_0, \rho )$ we set   $ z= T(y)$. Since $ | x_0| < 2r$   and 
$ | y-x_0| < r$, we estimate 
\begin{equation}
| z|= | T(y)|= \frac{r^2}{| y|}  \ge \frac{r^2}{ | y- x_0|+ | x_0|} >  \frac{r}{3}  .   
\label{8.12}
\end{equation}

 Making change of variable $z=T(y)$, using the transformation formula  of the Lebesgue integral  together with  \eqref{8.3},  \eqref{8.12} and \eqref{8.11}, we  infer 
\begin{align*}
II &= \frac{1}{c_n^2 \rho ^{ 2n}}\int_{B( x_0  , \rho)\cap B(r) }\int_{T(B( x_0 , \rho )  \setminus  B(r))} |u(z)-u  (y')| \Big(\frac{r}{| z|}\Big)^{ 2n} dz dy' 
\\
&\le c\intmw_{B(\hat{x}, 2\rho)\cap B(r) }\intmw_{B(\hat{x}, 2\rho )\cap B(r)} |u (z)-u(y')| dzdy' \le c \| u\|_{ BMO(B(r))}. 
\end{align*}
Similarly, we also  estimate
\[
III \le c \| u\|_{ BMO(B(r))}. 
\]
Therefore,  we have
\[
\intmw_{B(x_0, \rho )} | U  - U_{ B(x_0, \rho )} | dy \le c \| u\|_{ BMO(B(r))}. 
\]  
This completes the proof of the lemma.  \hfill \Beweisende 

\vspace{0.5cm}  

\section{Multiplicative inequalities}
\label{sec:-8}
\setcounter{secnum}{\value{section} \setcounter{equation}{0}
\renewcommand{\theequation}{\mbox{\arabic{secnum}.\arabic{equation}}}}

 \vspace{0.5cm}  
 In our discussion below we shall make use of the following multiplicative inequalities involving cut-off functions, where   $B(r)\subset \Bbb R^3$,  and  $u=(u_1, u_2, u_3)$ is a vector field on $\Bbb R^3$. 
 \begin{lem}
 \label{lem8.4}
 Let  $ \psi \in C^{\infty}_{\rm c}(B(r))$, $ 0< r< +\infty$, such that $ 0 \le \psi \le 1$ in $ B(r)$. For all $ u\in W^{1,\, q}(B(r))$, $ 2< q < +\infty,$ with $ \nabla \cdot u=0$ a.e. in $ B(r)$ and 
 for  all $ m \ge \frac{5q-6}{2q}$ it holds 
 \begin{align}
 \| \nabla  u \psi ^{ m}\|_{ q} &\le c  \| (\nabla\times  u) \psi^m \|_{ q} 
 +c \| \nabla \psi \|_{ \infty}^{ a} \| u \psi ^{ m- a}\|_{ 2},
 \label{8.13a}
 \\
 \| u \psi ^{ m-k}\|_{ q} &\le c \| u \psi ^{ m- ka}\|_{ 2}^{ \frac{1}{a}} 
 \|  (\nabla\times  u) \psi^m \|_{ q} ^{ 1- \frac{1}{a}}
 + c\| \nabla \psi \|_{ \infty}^{ a-1} \| u \psi ^{ m- ka}\|_{ 2},
 \label{8.13b}
 \end{align}
 where 
 \[
a= \frac{5q-6}{2q},\quad  1 \le k \le \frac{m}{a}.
\]
 \end{lem}

 {\bf Proof}: Recalling that $ \nabla \cdot u=0$, we calculate
 \[
- \Delta ( u\psi ^m) = \nabla \times \nabla \times (u \psi^m )- \nabla (u\cdot \nabla \psi^m). 
\]
 Applying $ \nabla $ to both sides of the above identity, we find that
 \[
- \Delta (\nabla ( u\psi^m) ) = \nabla \nabla \times (\nabla \times (u \psi^m ))- \nabla \nabla (u \cdot \nabla \psi^m). 
\]
 Using  Calder\'on-Zygmund's inequality, we see that 
 \begin{align}
 \| \nabla u \psi ^m\|_{q}  & \le \| \nabla (u \psi ^m)\|_{q} + c \| \nabla \psi \|_{ \infty}\| u \psi^{ m-1} \|_{q}
 \cr
 &\le c\| \nabla \times (\psi^m u)\|_{ q} 
  + c \| \nabla \psi \|_{ \infty}\| u \psi^{ m-1} \|_{q}
  \cr
  & \le 
c\| (\nabla \times   u) \psi ^m\|_{q} + c\| \nabla \psi \|_{ \infty}\| u \psi^{ m-1} \|_{ q}. 
\label{8.cz}
 \end{align}
 Estimating the last term on the right-hand side of the above estimate by means of  \eqref{B.1} with $ n=3$, and applying Young's inequality,   
 we get \eqref{8.13a}. The second inequality easily follows \eqref{B.1} together with  \eqref{8.13a}. 
 \hfill \Beweisende

 \vspace{0.5cm}  
Arguing as in the proof of Lemma\,\ref{lem8.4}, we get the following 

\begin{lem}
\label{lem8.4a}
Let    $ \psi \in C^{\infty}_c(B(r))$, $ 0< r< +\infty$, with $ 0 \le  \psi \le 1$. 
 For every  $ u \in W^{1,\, 1}(B(r))$ such that  $ \nabla \cdot u =0$  and $ \nabla \times  u\in BMO(B(r))$
it holds
\begin{equation}
\| \nabla  u\psi ^{ 6}\|_{ BMO} \le c \| (\nabla \times u) \psi^{ 5} \|_{ BMO} + 
c  \Big(r^{ \frac{3}{2}}  \| \nabla \psi \|^{4}_{ \infty}  +   \| \nabla \psi \|^{ \frac{5}{2}}_{ \infty}\Big)\| u \psi \|_{ 2}.
\label{8.14}
\end{equation}

\end{lem} 
 
 {\bf Proof}:  As we have seen in the proof of Lemma\,\ref{lem8.4},  it holds in $ \R^{3}$
 \[
- \Delta \nabla ( u\psi^{ 6} ) = \nabla \nabla \times (\nabla \times (u \psi^{ 6} ))- \nabla \nabla (u \cdot \nabla \psi^{ 6}). 
\]
 Using  Calder\'on-Zygmund's inequality\cite{ste}, we find that 
 \begin{align}
 \| \nabla  u \psi ^{ 6}\|_{ BMO}   &\le c\| \nabla \times (\psi^{ 6} u)\|_{ BMO} 
  + c\| u \cdot \nabla \psi^{ 6} \|_{ BMO} 
  \cr
  & \le 
c\| (\nabla \times   u) \psi ^{ 6}\|_{ BMO} + c\| \nabla \psi \|_{ \infty}\| u \psi^{5} \|_{ \infty}. 
\label{8.14a}
 \end{align}
 On the other hand, in view of  \eqref{B.2} with $ n=3$, $ q = 6$  and $ m= 5$ we get 
   \begin{equation}
\| u \psi ^{ 5}\|_{ \infty} \le c\| u \psi \|_{ 2}^{ \frac{1}{4 } }   \| \nabla  u \psi ^{ 5}\|^
{ \frac{3}{4}}_{ 6} + 
c\| \nabla \psi \|_{ \infty}^{ \frac{3}{2}} \| u \psi\|_{ 2}. 
\label{8.14b}
\end{equation}
We estimate the  first term  on the right-hand side of \eqref{8.14b}  by \eqref{8.13a} with $ q= 6$ and $ m=5$. This together with \eqref{B.3} 
 gives 
 \begin{align}
\| u \psi ^{ 5}\|_{ \infty} &\le c\| u \psi\|_{ 2}^{ \frac{1}{4} }   \| (\nabla \times  u) \psi ^{ 5}\|^
{ \frac{3}{4}}_{ 6} + c\| \nabla \psi \|_{ \infty}^{ \frac{3}{2}} \| u \psi \|_{ 2}. 
 \cr
 &\le c r^{ \frac{3}{8}}\| u \psi\|_{ 2}^{ \frac{1}{4 } }   \|( \nabla \times  u) \psi ^{ 5}\|^
{ \frac{3}{4}}_{ BMO(B(r))} + c\| \nabla \psi \|_{ \infty}^{ \frac{3}{2}} \| u \psi \|_{ 2}. 
 \label{8.14c}
 \end{align}
Finally, combining \eqref{8.14a} and \eqref{8.14c}, and applying Young's inequality, we obtain \eqref{8.14}. \hfill \Beweisende   

\vspace{0.3cm}

\begin{lem}
\label{lem8.4b} Let $ u \in W^{1,\, 1}(B(r))$ with $ \nabla \times u \in BMO(B(r))$. Then for all $ \psi \in C^{\infty}_{\rm c}(B(r))$ 
with $ 0 \le \psi \le 1$ we get 
\begin{align}
\| (\nabla \times u) \psi ^{5}\|_{ BMO} 
\le c \Big(1+ r^{ 5}\| \nabla \psi \|_{ \infty}^{ 5}\Big) \Big(| \nabla \times u|_{ BMO(B(r))} + cr^{- \frac{5}{2}} \| u\|_{ L^2(B(r))}\Big). 
\label{8.15zz} 
\end{align}

\end{lem}

{\bf Proof}:
Assume $ r=1$.  Let $ \eta \in C^{\infty}_{\rm c} (B(1))$ such that $ | \nabla \eta | \le c$ and $ \intl_{B(1)} \eta dx \ge c$, where 
$ c >0$ stands for a constant depending only on $ n$. For $ f\in L^1(B(1))$ we define the mean 
\[
\widetilde{f}_{ B(1)} = \frac{1}{ \intl_{B(1)} \eta dx } \intl_{B(1)} f\eta dx.  
\]
First,  using triangle inequality and applying integration by parts,  we see that 
\begin{align}
\| \nabla \times u\|_{ L^1(B(1))} &= \| \nabla \times u -\widetilde{ \nabla \times u}_{ B(1)}  \|_{ L^1(B(1))} + 
| \widetilde{ \nabla \times u}_{ B(1)}|
\cr
&\le c| \nabla \times u|_{ BMO(B(1))} + c \| u\|_{ L^1(B(1))}.
\label{l1}
\end{align}

Using \eqref{l1}, we  estimate  for $ \rho \ge \frac{1}{2}$ and $ x_0\in \R^{3}$
\begin{align*}
& \intmw_{B(x_0, \rho )} 
| (\nabla \times u) \psi^5- ((\nabla \times u) \psi^5)_{ B(x_0, \rho )}|  dx
\\
&\qquad \le c  \| \nabla \times u\|_{ L^1(B(1))} 
\\
&\qquad \le 
c| \nabla \times u|_{ BMO(B(1))} + c \| u\|_{ L^1(B(1))}.
\end{align*}
In case $ \rho \le \frac{1}{2}$ and $ B(x_0, \rho ) \cap B(1) \neq \emptyset$ 
there exists $ y_0 \in B(1)$ such that $ B(x_0, \rho ) \subset B(y_0, 2\rho )$ and 
\begin{align}
& \intmw_{B(x_0, \rho )} | (\nabla \times u) \psi^5 - 
((\nabla \times u) \psi^5)_{ B(x_0, \rho )}|  dx
\cr
&\le c\intmw_{B(y_0, 2\rho )} \intmw_{B(y_0, 2\rho )} 
| (\nabla \times u)(x) \psi^5(x)- |  (\nabla \times  u)(y) \psi^5(y)
|  dxdy
\cr
&\le c| \nabla \times u|_{ BMO(B(1))} + c \intmw_{B(y_0, 2\rho )} \intmw_{B(y_0, 2\rho )} | (\nabla \times u)(x)  
|  | \psi^5(x) - \psi^5(y)|dxdy.
 \label{8.bmo1}
\end{align}
By the fundamental theorem of differentiation  and integration  we find 
$ \xi _i \in [x,y]$, $ i=1, 2, 3$ such that with $\xi _{0} =y$
\begin{align*}
\psi ^5(y) -\psi ^5(x) &\le 5  \psi ^{ 4}(\xi_1 )\nabla \psi (\xi_1 )\cdot (y-x) 
\\
&= 5\psi ^{ 4}(x)\nabla \psi (\xi _1)\cdot (y-x) +5  (\psi ^{ 4}(\xi_1 )-\psi ^{ 4}(x))\nabla \psi (\xi_1 )\cdot (y-x) 
\\
&= \sum_{k=1}^{2}\psi ^{ 5-k}(x) \prod_{ i=1}^k (6-i) \nabla \psi (\xi _i)\cdot (\xi _{ i-1}-x)
\\
&\qquad + \psi ^{ 2}(\xi _3) \prod_{ i=1}^3 (6-i) \nabla \psi (\xi _i)\cdot (\xi _{ i-1}-x).
\end{align*}
  This   yields 
\begin{align*}
&\intmw_{B(y_0, 2\rho )} \intmw_{B(y_0, 2\rho )} | \nabla \times u(x)  
|  | \psi^5(x) - \psi^5(y)|dxdy
\\
&\qquad \le c\sum_{k=1}^{2}  \| \nabla \psi \|^k_{ \infty}\rho^{ k-3}  \intl_{B(y_0, 2\rho )} | \nabla \times u | \psi^{ 5-k} dx 
 + c \| \nabla \psi \|_{ \infty}^3 \| \nabla \times u\|_{ L^1(B(1)).}
\end{align*}
Taking into account Lemma\,\ref{lemB.3} we have $\nabla \times u\in L^{q}(B(1))$ for all $1 \le q <+\infty $. Applying H\"older's inequality, we obtain
\begin{align}
&\intmw_{B(y_0, 2\rho )} \intmw_{B(y_0, 2\rho )} | \nabla \times u(x)  
|  | \psi^5(x) - \psi^5(y)|dxdy
\cr
&\qquad \le c\sum_{k=1}^{2}  \| \nabla \psi \|^k_{ \infty} \| (\nabla \times u) \psi ^{ 5-k}\|_{ \frac{3}{k}} 
+ c \| \nabla \psi \|_{ \infty}^3 \| \nabla \times u\|_{ L^1(B(1))}. 
 \label{8.bmo2}
\end{align}
Esimating the right-hand side of  \eqref{8.bmo1} by means of  \eqref{8.bmo2}, we already verify that  that $(\nabla \times u)\psi ^{5}\in BMO$. By using H\"older's inequality, together  with John-Nirenberg's inequality, we estimate 
\begin{align*}
\| (\nabla \times u) \psi ^{ 5-k}\|_{ \frac{3}{k}}  &\le \|\nabla \times u\|_{ L^{1}(B(1))}^{ \frac{2k-1}{5}}  
\| (\nabla \times u) \psi ^{ (5-k) \frac{5}{6-2k}}\|_{ 6}^{ \frac{6-2k}{5}}
\\
&\le \|\nabla \times u\|_{ L^{1}(B(1))}^{ \frac{2k-1}{5}}  
\| (\nabla \times u) \psi ^{5}\|_{6}^{ \frac{6-2k}{5}}
\\
&\le c \|\nabla \times u\|_{ L^{1}(B(1))}^{ \frac{2k-1}{5}}  
\| (\nabla \times u) \psi ^{5}\|_{ BMO}^{ \frac{6-2k}{5}}.
\end{align*}

Combining the above inequalities, and applying Young's inequality together with \eqref{l1}, we arrive at 
\begin{align}
&\| (\nabla \times u) \psi ^{5}\|_{ BMO} 
\cr
& \le c | \nabla \times u|_{ BMO(B(1))}+ c \sum_{k=1}^{2} \| \nabla \psi \|_{ \infty}^{ \frac{5k}{2k-1}} \| \nabla \times u\|_{ L^1(B(1))} 
+ c \| \nabla \psi \|_{ \infty}^3 \| \nabla \times u\|_{ L^1(B(1))}
\cr
&\le c (1+ \| \nabla \psi \|_{ \infty}^{ 5}) \Big(| \nabla \times u|_{ BMO(B(1))} + c\| u\|_{ L^1(B(1))}\Big). 
\label{8.15z} 
\end{align}
Whence, \eqref{8.15zz} follows immediately from \eqref{8.15z} by using a standard scaling argument.  

\hfill \Beweisende

Combining Lemma\,\ref{lem8.4a} and Lemma\,\ref{lem8.4b}, we get 

\begin{cor}
\label{cor8.5} 
Let    $ \psi \in C^{\infty}_c(B(r))$, $ 0< r< +\infty$, with $ 0 \le  \psi \le 1$.
 For every  $ u \in W^{1,\, 1}(B(r))$ such that  $ \nabla \cdot u =0$  and $ \nabla \times  u\in BMO(B(r))$
it holds
\begin{align}
\| \nabla u \psi ^{ 6}\|_{ BMO(B(r))} 
&\le c\Big\{1+ r^{ 5}\| \nabla \psi \|_{ \infty}^{ 5}\Big\}| \nabla \times u|_{ BMO(B(r))} 
\cr
&  \qquad +c \Big\{r^{ -\frac{5}{2}} + r^{ \frac{5}{2}}\| \nabla \psi \|^{ 5}_{ \infty}\Big\}\| u\|_{ L^2(B(r))}.
\label{8.15e}
\end{align}
\end{cor} 

{\bf Proof}: Combining \eqref{8.14} and \eqref{8.15zz} along with Young's inequality, we infer 
\begin{align}
&\| \nabla u \psi ^{ 6}\|_{ BMO(B(r))} 
\cr
&\le c \Big\{1+ r^{ 5}\| \nabla \psi \|_{ \infty}^{ 5}\Big\} (| \nabla \times u|_{ BMO(B(r))} + cr^{ - \frac{5}{2}} \| u\|_{ L^2(B(r))})
\cr
&  \qquad\qquad  +c  \Big\{r^{ \frac{3}{2}} \| \nabla \psi \|^{ 5}_{ \infty}+   \| \nabla \psi \|^{ \frac{5}{2}}_{ \infty}\Big\}\| u \|_{ L^2(B(r))}
\cr 
&  \le c \Big\{1+ r^{ 5}\| \nabla \psi \|_{ \infty}^{ 5}\Big\}| \nabla \times u|_{ BMO(B(r))} 
\cr
&\qquad\qquad   + c \Big\{r^{ -\frac{5}{2}} + r^{ \frac{5}{2}}\| \nabla \psi \|^{ 5}\Big\}\| u\|_{ L^2(B(r))}.
\label{8.15l}
\end{align}
Whence, \eqref{8.15e}.  \hfill \Beweisende  

\begin{rem}
\label{rem8.9}
Note that thanks to \eqref{8.15e} for every   $ v\in L^\infty(-\rho , 0; L^2(B(\rho )))$,
which satisfies  the  local Beale-Kato-Majda condition  \eqref{bkm}  it holds for all $ 0< r < \rho $
\begin{equation}
 \intl_{-1}^{0} \| \omega(t) \|_{ BMO(B(r))}  dt <+\infty. 
\label{bkm1}
\end{equation}
In particular, \eqref{8.14c} with $ n=3$ together with \eqref{bkm1}  implies that for all $ 0< r< \rho $
\begin{equation}
\label{vel}
v\in L^{ \frac{4}{3}} (-\rho , 0; L^\infty(B(r))). 
\end{equation}

\end{rem}

\begin{lem}
\label{lem8.6}
Let $ u \in W^{2,\, q}(B(r ))$, $ 2 \le q < +\infty$. Let $ m,k\in \R$ such that  $2 \le   m <+\infty$ and $ 0 < k \le 2m$. 
Then for every $ \psi \in C^{\infty}_c(B(r))$ it holds  
\begin{equation}
\| (\nabla u) \psi^m \|_{ q} 
 \le c\| u \psi ^{ 2m-k} \|_{ q}^{ 1/2} \|(\nabla^2  u)\psi ^{k}\|_{ q}^{1/2} +
  c \| \nabla \psi \|_{ \infty}  \| u \psi ^{ m-1}\|_{ q}. 
\label{8.16}
\end{equation}
If in addition, if $ \nabla \cdot u=0$ almost everywhere in $ B(r )$, then for $2  \le  m < +\infty$ and $ m+1 \le k \le 2m$
it holds
\begin{align}
\|( \nabla^2 u) \psi^k \|_{ q} 
& \le c \|(\nabla \nabla \times  u) \psi ^{k}\|_{ q} + c \| \nabla \psi \|^2_{ \infty} 
 \| u \psi ^{ k-2}\|_{ q}, 
\label{8.17}
\\
\| (\nabla u) \psi^m \|_{ q} 
 &\le c\| u  \psi^{ 2m-k}\|_{ q}^{ 1/2} \|(\nabla \nabla\times   u) \psi ^{k}\|_{ q}^{1/2} + 
 c \| \nabla \psi \|_{ \infty}  \| u\psi ^{ m-1}\|_{ q}. 
\label{8.17a}
\end{align}

\end{lem}

{\bf Proof}:  
Applying integration by parts, and using H\"older's inequality, we get 
\begin{align*}
\|(\nabla  u) \psi^m \|^{ q}_{q}  
& =-\intl_{B(\rho )} u \nabla u\cdot  (\nabla  | \nabla u|^{ q-2}) \psi ^{ mq} dx- 
\intl_{B(\rho )} u \nabla u  | \nabla u|^{ q-2} \cdot \nabla \psi ^{ mq} dx
\\
& \le c \| u \psi ^{ 2m-k}\|_{q}  \| (\nabla ^2 u )\psi ^{ k}\|_{q}  \|(\nabla  u) \psi^m \|^{ q-2}_{q}  
\\
&\qquad + c \| \nabla \psi \|_{ \infty} \| u \psi^{ m-1}\|_{L^q(B(\rho ))}   
\|(\nabla  u) \psi^{ m} \|^{ q-1}_{q},  
\end{align*}  
and Young's inequality gives \eqref{8.16}.   

\hspace{0.5cm}
In case $ \nabla \cdot u =0$ almost everywhere in $ B(r )$ we may apply \eqref{8.cz} with $ \nabla u$ in place of $ u$  in order to estimate in \eqref{8.16} 
the norm  involving the second gradient of $ u$. This together with \eqref{8.16} with $ m=k-1$  gives 
\begin{align*}
\| (\nabla ^2 u )\psi^k \|_q &\le c \| (\nabla  \nabla \times u )\psi^k \|_q + 
c\| \nabla \psi \|_{ \infty}\| \nabla u\psi^{ k-1} \|_{ q} 
\\
&\le  c \| (\nabla  \nabla \times u) \psi^k\|_q + c\| \nabla \psi \|_{ \infty}\| u \psi^{ k-2} \|^{ 1/2}_{q} 
\|( \nabla^2 u)\psi ^{ k}\|^{ 1/2}_{ q}  
\\
&\qquad \qquad + \| \nabla \psi \|_{ \infty}^2\| u \psi^{ k-2}\|_{ q}. 
\end{align*}
Then we  apply Young's inequality to obtain \eqref{8.17}.  The estimate \eqref{8.17a} is now an immediate consequence of \eqref{8.16} 
and \eqref{8.17}.   
\hfill \Beweisende  

\vspace{0.3cm}
Combining Lemma\,\ref{lem8.4} and Lemma\,\ref{lem8.6}, we get the following  

\begin{cor}
\label{cor8.7} 
For all $ u \in W^{2,\, q}(B(r))$, $ 2 < q < +\infty $,  for all $ \psi \in C^{\infty}_{\rm c}(B(r))$ with $0 \le \psi  \le 1$ and for all $ m >5$ we get 
 \begin{equation}
\| (\nabla ^2 u) \psi ^{m}\|_{ q} \le c \| (\nabla \nabla \times u) \psi ^m\|_{ q} +  c\| \nabla \psi \|_{ \infty}^{a+1} \| u \psi ^{ m -5}\|_{2},
\label{8.22}
\end{equation} 
 where
\[
a= \frac{5q- 6}{2q}. 
\]
\end{cor}

{\bf Proof}:  Let $ m >5$. Using  \eqref{8.17} with $k=m$, we get 
\begin{align}
\|( \nabla^2 u) \psi^m \|_{ q} 
& \le c \|(\nabla \nabla \times  u) \psi ^{m}\|_{ q} + c \| \nabla \psi \|^2_{ \infty} 
 \| u \psi ^{ m-2}\|_{ q}, 
\label{8.17c}
\end{align}

  Next, in view of  \eqref{8.13b}  with $ k=2$ and  \eqref{8.17a} with $k=m$, we get 
 \begin{align*}
\| u \psi ^{ m-2}\|_{ q} &\le c \| u \psi ^{ m- 2a}\|_{ 2}^{ \frac{1}{a}} 
\| (\nabla\times  u) \psi^{ m} \|_{ q} ^{ 1- \frac{1}{a}}
 + c\| \nabla \psi \|_{ \infty}^{ a-1} \| u \psi ^{ m- 2a}\|_{ 2},
 \\
 &\le c \| u \psi ^{ m- 2a}\|_{ 2}^{ \frac{1}{a}} 
 \| u \psi ^{m }\|^{ \frac{1}{2}- \frac{1}{2a}}_{q} 
 \|(\nabla \nabla \times u) \psi ^{m}\|^{ \frac{1}{2}- \frac{1}{2a}}_{q}  
 \\
&  \qquad + c \|\nabla \psi \|_{\infty }^{ 1- \frac{1}{a}} \| u \psi ^{ m- 2a}\|_{ 2}^{ \frac{1}{a}}  \| u \psi ^{m -1}\|^{ 1- \frac{1}{a}}_{q} 
 \\
  &\qquad + c\| \nabla \psi \|_{ \infty}^{ a-1} \| u \psi ^{ m- 2a}\|_{ 2}.
 \end{align*} 

Applying Young's inequality, noting that $ \psi ^{m-2a} \le \psi ^{m-5}$ due to $2a \le 5$,  and combining the resultant inequality with  \eqref{8.17c}, 
we obtain \eqref{8.22}.  \hfill \Beweisende

\section{ Local estimates of the vorticity }
\label{sec:-?}
\setcounter{secnum}{\value{section} \setcounter{equation}{0}
\renewcommand{\theequation}{\mbox{\arabic{secnum}.\arabic{equation}}}}

\begin{thm}
\label{thm10.6}
Under the assumption of Theorem\,\ref{thm1.1} it holds for all $ 0< r < \rho $ and for all $ 1 \le q <+\infty$
\begin{align}
\omega \in L^\infty(-\rho , 0; L^q(B(r))).
\label{10.41}
\end{align}

\end{thm}

{\bf Proof}:  Applying $ \curl$ to both sides of  \eqref{euler}, we get the vorticity equation 
\begin{equation}
\partial _t \omega + v\cdot \nabla \omega = \omega \cdot \nabla v\quad  \text{ in}\quad  B(\rho )\times (-\rho , 0). 
\label{vorteq}
\end{equation}
Fix, $ 2 \le q < +\infty$. 
Let $ 0< r < \sigma < \rho $ be arbitrarily chosen. Let $ 0< r< r_1 < r_2 < \sigma $, and set $ \widetilde{r} := \frac{r_1+r_2}{2}$.
Let $ \phi\in C^{\infty}_{\rm c}(B(\widetilde{r} )) $ 
denote a cut off function such that $ 0 \le \phi  \le 1$,  $ \phi \equiv 1$ on $ B(r) $, and $ | \nabla \phi | \le c (r_2 -r_1)^{ -1}$. 
Next, let $  -\rho < t_0 <t_1 < 0$ be fixed, where $ t_0 $ will be taken sufficiently small, specified below. 
We multiply \eqref{vorteq} by $ \omega | \omega |^{ q-2} \phi^{ q}$, integrate the result over $ B(r_2)\times (t_0, t )$, 
$ t_0 < t < t_1$, and apply integration by parts. This yields  
\begin{align}
&\| \omega (t) \phi \|_{ q}^q 
\cr
&= \| \omega (t_0) \phi \|_{ q}^q +
 \intl_{t_0 }^{t} \intl_{B(r_2)}    v(s)\cdot \nabla \phi  | \omega(s) |^q \phi ^{ q -1} dx ds 
 \cr
&\qquad \qquad  + 
  q\intl_{t_0}^{t}  \intl_{B(r_2)}  \omega(s) \cdot \nabla v(s) \cdot  \omega(s) | \omega(s) |^{ q-2} \phi ^{ q } dx ds
 \cr
 & = \| \omega (t_0) \phi\|_{ q}^q + I+ II. 
 \label{10.48}
\end{align}
 First, by using H\"older's inequality and Lemma \ref{lemB.3}, we estimate 
 \begin{align*}
I &\le c (r_2 -r_1)^{ -1}  \intl_{t_0}^{t} \| v(s)\|_{ L^\infty(B(\sigma ))} \| \omega(s) \phi^{ 1 - \frac{1}{q}} \|^q_{q } ds
 \\
 & \le  c (r_2 -r_1)^{ -1}  \intl_{t_0}^{t} \| v(s)\|_{ L^\infty(B(\sigma ))} \| \omega(s) \phi^{ 1 - \frac{1}{q}} \|^{ \frac{1}{4}}_{q } ds
\esssup_{s\in (t_0, t_1)} \| \omega(s) \|^{ q- \frac{1}{4}}_{ L^q(B(r_2))}
\\
 & \le  c (r_2 -r_1)^{ -1} \| v\|_{ L^{ \frac{4}{3}}(-\rho , 0; L^\infty(B(\sigma )))}  
 \bigg(\intl_{-\rho }^{0} \| \omega (s)\|_{ BMO(B(\sigma ))}ds\bigg)^{ \frac{1}{4}}\times 
\\
&\qquad \times \esssup_{s\in (t_0, t_1)} \| \omega(s) \|^{ q- \frac{1}{4}}_{ L^q(B(r_2))}. 
 \end{align*}
 Using Young's inequality, we obtain 
 \begin{align*}
 I &\le c (r_2 -r_1)^{ -4q} \| v\|^{ 4q}_{ L^{\frac43}(-\rho , 0; L^\infty(B(\sigma )))} 
 \bigg(\intl_{-\rho }^{0} \| \omega (s)\|_{ BMO(B(\sigma  ))}  ds\bigg)^{ q} 
\\
&\qquad +\var \esssup_{s\in (t_0, t_1)} \| \omega(s) \|^{ q}_{ L^q(B(r_2))}. 
 \end{align*}
Secondly, we get by the aid of H\"older's inequality 
\begin{align}
II \le c  \intl_{t_0}^{t}  \| \omega(s) \cdot \nabla v(s) \|_{L^ q(B(\widetilde{r} ))}    ds
\esssup_{s\in (t_0, t_1)} \| \omega(s) \|^{ q-1}_{ L^q(B(r_2))}.
\label{10.49}
\end{align}

To proceed further we prove the following localization of Kozono-Taniuchi's inequality\cite{koz1}.

\begin{lem}
Let $ f, g \in BMO(B(r))\cap L^q(B(r))$, $ 1<q<+\infty$, then $ f\cdot  g \in L^q(B(r))$ and it holds 
\begin{align}
\| f\cdot g\|_{ L^q(B(r))} &\le 
c\Big(| f|_{ BMO(B(r))} \| g\|_{ L^q(B(r))} + | g|_{ BMO(B(r))} \| f\|_{ L^q(B(r))}\Big)
\cr
& \qquad +c r^{- \frac{3}{q}}\| f \|_{L^q(B(r))} \| g\|_{ L^q(B(r))},
\label{10.40a}
\end{align}
where the constant $ c>0$ depends on $ q$ only. 
\end{lem}

{\bf Proof}: We define the extensions $\widetilde{f}$ and $ \widetilde{g}$ of $f$ and $g$ respectively as follows. 
\[
\widetilde{f}(x) = 
\begin{cases}
f(x)\quad  \text{ if} \quad  x\in B(r)
\\[0.3cm]
f \Big(\frac{r^2 x}{ | x|^2}\Big) \phi (x)\text{ if} \quad  x\in B(r)^c. 
\end{cases} 
\]
\[
\widetilde{g}(x) = 
\begin{cases}
g(x)\quad  \text{ if} \quad  x\in B(r)
\\[0.3cm]
g \Big(\frac{r^2x}{ | x|^2}\Big) \phi (x)\text{ if} \quad  x\in B(r)^c, 
\end{cases} 
\]
where $ \phi \in C^{\infty}_{\rm c}(B(4r))$ denotes a suitable cut off function such that $ \phi \equiv 1$ on $ B(2r)$. 
According to Lemma\,\ref{lem8.1} it holds $ f, g \in BMO \cap L^q(\R^{3})$ together with the estimates 
\begin{align}
&\| \widetilde{f} \|_{ q} \le c \| f\|_{ L^q(B(r))},\quad  \| \widetilde{g} \|_{ q} \le c \| g\|_{ L^q(B(r))}
\label{10.40b}
\\
&\quad \| \widetilde{f} \|_{ BMO} \le c | f|_{ BMO(B(r))} + c r^{ -3} \| f\|_{ 1},\quad  
\label{10.40c}
\\
&\quad \| \widetilde{g} \|_{ BMO}\le c | g|_{ BMO(B(r))} + c r^{ -3} \| g\|_{ 1}. 
\label{10.40d}
\end{align}

Thanks to Kozono-Taniuchi's inequality\cite[Lemma 1(i)]{koz1},  combined  with \eqref{10.40b}, \eqref{10.40c} and \eqref{10.40d} we find that 
\begin{align*}
\| f\cdot g\|_{ L^q(B(r))} &\le \| \widetilde{f} \cdot \widetilde{g} \|_{ q} 
\\
&\le c\Big(\| \widetilde{f} \|_{ BMO} \| \widetilde{g} \|_{q} + \| \widetilde{g} \|_{ BMO} \| \widetilde{f} \|_{ q}\Big)
\\
&\le c\Big(| f |_{ BMO(B(r))} \| g\|_{ L^q(B(r))} + |g |_{ BMO(B(r))} \| f\|_{ L^q(B(r))}\Big)
\\
&\qquad \qquad +c r^{ -3}\Big( \| f \|_{L^1(B(r))} \| g\|_{ L^q(B(r))} + \|g \|_{ L^1(B(r))} \| f\|_{ L^q(B(r))}\Big).
\end{align*} 
Whence, by using Jensen's inequality we get  \eqref{10.40a}.  \hfill \Beweisende  

\vspace{0.5cm}  
Thanks to \eqref{10.40a} we find   
\begin{align}
&\| \omega(s) \cdot \nabla v(s)\|_{ L^q(B(\widetilde{r} ))} 
\cr
&\qquad \le 
c\Big(| \omega(s)| _{ BMO(B(\widetilde{r}))} \| \nabla v(s)\|_{ L^q(B(\widetilde{r}))} + |  \nabla v(s)|_{ BMO(B(\widetilde{r}))} \| \omega (s)\|_{ L^q(B(\widetilde{r}))}\Big)
\cr
& \qquad \qquad +c r^{-\frac{3}{q}} \| \omega (s) \|_{L^q(B(\widetilde{r}))} \| \nabla v(s)\|_{ L^q(B(\widetilde{r}))}.
\label{10.40e}
\end{align}
By the aid of Corollary\,\ref{cor8.5} we get 
\begin{equation}
 \| \nabla v(s)\|_{ BMO(B(\widetilde{r}))} \le c\Big(| \omega  (s)|_{ BMO(B(\rho ))}  + (r_2-r_1)^{-\frac52}\| v(s)\|_{ L^2(\rho )}\Big). 
\label{10.40f}
\end{equation}
Furthermore, by virtue of Lemma\,\ref{lem8.4} we see that 
\begin{equation}
\| \nabla v(s)\|_{ L^q(B(\widetilde{r}))} \le c\Big(\| \omega  (s)\|_{L^q(B(r_2 ))}  + 
(r_2-r_1)^{  -\beta  }\| v(s)\|_{ L^2(\rho )}\Big),
\label{10.40g}
\end{equation}
where
\[
\beta =  \frac{5q-6}{2q}.
\]
We now estimate the right-hand side of \eqref{10.40e} by \eqref{10.40f} and \eqref{10.40g}. This yields 
\begin{align*}
&\| \omega(s) \cdot \nabla v(s)\|_{ L^q(B(\widetilde{r} ))} 
\\
&\quad   \le c |\omega (s)|_{ BMO(B(\rho ))} \| \omega (s)\|_{ L^q(B(r_2))}
\\
&\qquad +c(r_2-r_1)^{-\frac52}\| v (s)\|_{ L^2(B(\rho ))} \| \omega (s)\|_{ L^q(B(r_2))} \\
 &\qquad+ c (r_2-r_1)^{ -\beta }| \omega (s)|_{ BMO(B(\rho ))} \| v(s)\|_{ L^2(B(\rho ))}
+ c (r_2-r_1)^{ -\beta -\frac52}\| v(s)\|_{ L^2(B(\rho ))}^2.
\end{align*}   
Noting that 
\[
\| \omega (s)\|_{ L^q(B(r_2))} \le c(r_2-r_1)^{\frac{3}{q}}\Big(| \omega (s)|_{ BMO(B(\rho ))} + (r_2-r_1)^{-\frac52}\| v(s)\|_{ L^2(B(\rho ))}\Big)
\]
due to \eqref{lemB.3},
we are led to 
\begin{align}
&\| \omega(s) \cdot \nabla v(s)\|_{ L^q(B(\widetilde{r} ))} 
\cr
& \qquad \le c a(s) \| \omega (s)\|_{ L^q(B(r_2))} + c (r_2- r_1)^{ -\beta  } a(s) \| v(s)\|_{ L^2(B(\rho ))}, 
\label{10.50}
\end{align}
where we have  set
\[
a(s) = | \omega (s)|_{ BMO(B(\rho ))}  +(r_2-r_1)^{-\frac52} \| v (s)\|_{ L^2(B(\rho ))}.
\]
Integrating \eqref{10.50} over $ (t_0, t)$, then combining the result with \eqref{10.49}, and applying Young's inequality, we obtain   
\begin{align*}
II \le c  \intl_{t_0}^{0} a(s)   ds \esssup_{s\in (t_0, t_1)} \| \omega(s) \|^{ q}_{ L^q(B(r_2))}
+ \var \esssup_{s\in (t_0, t_1)} \| \omega(s) \|^{ q}_{ L^q(B(r_2))} 
\\
+ c (r_2-r_1)^{ -q\beta }\| v\|_{ L^\infty(-\rho , 0; L^2(B(\rho )))}^q \bigg(\intl_{t_0}^{0} a(s)   ds\bigg)^q. 
\end{align*}
We take $ \var = \frac{1}{6}$, and choose $ t_0$ such that 
\[
c \intl_{t_0}^{0} a(s)  ds  \le \frac{1}{6}.
\]
Inserting the estimates of $ I$ and $ II$ into the right-hand side of \eqref{10.48}, and taking into account that $ \beta  \le 4$,  we get 
\begin{align}
&\esssup_{s\in (t_0, t_1)} \| \omega(s) \|^{ q}_{ L^q(B(r_1))} 
\cr
 & \quad \le \frac{1}{2}\esssup_{s\in (t_0, t_1)} \| \omega(s) \|^{ q}_{ L^q(B(r_2))} 
\cr
&\qquad +\| \omega (t_0) \|_{ L^q(B(\sigma ))}+ c (r_2-r_1)^{ -4q }\Big(\| v\|_{ L^\infty(-\rho , 0; L^2(B(\rho )))}^q + \| v\|_{ L^{\frac43}(-\rho , 0; L^\infty(B(\rho )))}^{ 4q}\Big).
\label{10.51}
\end{align}
By using a standard iteration argument, we deduce from \eqref{10.51} that 
\begin{align*}
&\esssup_{s\in (t_0, t_1)} \| \omega(s) \|^{ q}_{ L^q(B(r))} 
\\
&\qquad \le c (\sigma -r)^{ -4q }\Big(\| \omega (t_0) \|_{ L^q(B(\sigma ))}+\| v\|_{ L^\infty(-\rho , 0; L^2(B(\rho )))}^q +
 \| v\|_{ L^{\frac43}(-\rho , 0; L^\infty(B(\rho )))}^{ 2q}\Big).  
\end{align*}
 Since $ c$ is independent of $ t_1$, using the fact $ \omega \in L^\infty(-\rho , t_0; L^q(B(r)))$ and \eqref{vel}, combined with the hypothesis
 $v\in L^\infty(-\rho, 0; L^2 (B(\rho)))$, 
 we get \eqref{10.41}.  \hfill \Beweisende 
 
 \vspace{0.3cm}
 By means of Sobolev's embedding theorem we immediately deduce from \eqref{10.41} the following 
 
 \begin{cor}
 \label{cor10.1}
 Let $ v\in L^\infty(-\rho , 0; L^2(B(\rho )))$ with $ |\omega(\cdot )|_{ BMO(B(\rho ))}\in L^1(-\rho ,0)$ be a solution to the 
 Euler eqautions \eqref{euler}. Then for all $ 0<r<\rho $, 
 \begin{equation}
 v\in L^\infty(-\rho, 0; C^{ 0, \gamma }(B(r)))\quad  \forall\,0< \gamma <1.
 \label{10.52}
 \end{equation} 
 In particular, $ v\in L^\infty(B(r)\times (-\rho , 0))$ for all $ 0<r<\rho $. 
 
 \end{cor}

 \section{Proof of Theorem \ref{thm1.1} }
 \label{sec:-?}
 \setcounter{secnum}{\value{section} \setcounter{equation}{0}
 \renewcommand{\theequation}{\mbox{\arabic{secnum}.\arabic{equation}}}}
 We are now ready to prove our main theorem.\\
 
 \noindent{\bf Proof of Theorem \ref{thm1.1}}:   We take  $-\rho < t_{ \ast} <0$ such that 
 \begin{equation}
 \frac{1}{2} \rho ^2  (-t_{ \ast})^{- \frac{1}{2}  }- 4C_0>0. 
 \label{tast}
 \end{equation}
 Let $ t_{ \ast}< t_0 < 0$. We consider the following transformation of $(v,p) \mapsto (V,P)$ defined by
 \begin{align*}
V (y,t) &= v((1+ (-t)^{  \frac{1}{2} } ) y, t),
\\
P (y,t) &= \frac{1}{1+(-t)^{  \frac{1}{2}  }}p((1+ (-t)^{  \frac{1}{2} }  ) y, t),
\end{align*}
which was  first introduced by the authors of this paper in \cite{cha1}.
Thanks to   Corollary\,\ref{cor10.1} we have 
\begin{equation}
\| V(t)\cdot y \|_{ L^\infty(B(\rho _0) )} \le C_0 \quad \forall\,t\in [t_0,0),
\label{7.a}
\end{equation} 
where $ \rho _0 := \frac{\rho }{(1+ (-t_0)^{  \frac{1}{2}  })} > \frac{\rho }{2}$.
We also define  
\[
W(y,t) =  \frac{  \frac{1}{2}   (-t)^{ -1/2}y + V(y,t)}{1+ (-t)^{ 1/2}},\quad  (y,t)\in B(\rho _0)\times (t_0, 0). 
\]
We claim that   
\begin{equation}
 W(y,t)\cdot y >0 $ for all $ (y,t)\in  \left(\overline{B (\rho _0)}  \setminus B(\rho /2)\right)\times (t_0, 0).\label{7.a2}
\end{equation}
Indeed, according to \eqref{7.a} together with \eqref{tast}  we estimate   
\begin{align*}
 W(y,t)\cdot y &=
 \frac{1}{2}   (-t)^{ - \frac{1}{2}}| y|^2 + V(y, t)\cdot y \ge  \frac{1}{2}  (-t)^{-  \frac{1}{2} } \frac{  \rho^2}{4} + V(y, t)\cdot y
\\
&\ge \frac{1}{4} \Big(  \frac{1}{2} \rho ^2  (-t)^{ -  \frac{1}{2} }  - 4C_0 \Big)\ge \frac{1  }{4} \Big(  \frac{1}{2}  \rho ^2 (-t_{ \ast})^{   -  \frac{1}{2} }  - 4C_0\Big) >0.
\end{align*}

Using the chain rule, we see that \eqref{euler} turns into the following equations,  which hold in $ B(\rho _0)\times (t _0,0).$
  \begin{align}
\partial _t V  + W \cdot \nabla V &= -\nabla P, \quad  \nabla \cdot V=0.
\label{euler1}  
  \end{align}
  We set
  \[
\Omega = \nabla \times V\quad \text{ in}\quad    B(\rho _0 )\times (t_0, 0).
\]
Then  applying $ \nabla \times $ to \eqref{euler1},  we obtain the following equations   
   \begin{align}
 \partial _t \Omega  + \frac{  \frac{1}{2} (-t)^{ -  \frac{1}{2} }}{1+ (-t)^{  \frac{1}{2}  }} \Omega  + 
 W\cdot \nabla \Omega  &=\frac{1}{1+(-t)^{\frac12}} \Omega \cdot \nabla V,
  \label{VorEqu}
 \end{align}

Applying the operator $ \partial _i  (i\in \{ 1,2,3\})$ to both sides of \eqref{VorEqu}, we get the equations  
   \begin{align}
 \partial _t \partial _i \Omega  + \frac{ \frac{1}{2} (-t)^{ -  \frac{1}{2} }}{1+ (-t)^{  \frac{1}{2}  }} \partial _i\Omega  + 
 W\cdot \nabla \partial _i\Omega  &=\frac{1}{1+(-t)^{\frac12}}  (\partial _i\Omega) \cdot \nabla V+\frac{1}{1+(-t)^{\frac12}}  \Omega \cdot \partial _i \nabla V.
  \label{VorEqu1}
 \end{align}

Let $ 3< q<+\infty$ and $ \frac{\rho}{2}< r < \rho _0$ be arbitrarily chosen. Set $ \rho _{ \ast} = \frac{r+ \rho _0}{2}$, and define
\[
r_m := \rho _{ \ast}- (\rho_{ \ast}-r)^{ m+1} \rho _{ \ast}^{ -m}, \quad  m \in \N\cup \{ 0\}. 
\]
Clearly, 
\[
r_{m+1}- r_{ m} = r \left(1- \frac{r}{\rho_* }\right)^{ m+1},\quad\text{and}\quad r_m\nearrow \rho_*.
\]

Let $ \eta_m \in  C^{\infty} (\R)$ denote a cut off function such that $ 0 \le \eta_m  \le  1 $ 
in $ \R$, $ \eta_m \equiv  1$ on $ (-\infty, r_{ m}]$, $ \eta_m \equiv 0$ in $ (r_{ m+1}, +\infty)$, and 
$ 0 \le  -\eta '_m \le \frac{2}{r_{ m+1}- r_{ m}} = 2 r^{ -1} \Big(\frac{\rho _0+r}{\rho _0-r}\Big)^{m+1}$. 
 Next,  we multiply both sides of \eqref{VorEqu1} by $ \partial _i \Omega | \nabla \Omega |^{ q-2} \phi_m^{ 6q}$,  where 
 $\phi_m (y) = \eta_m (| y|) $, integrate the result over  $ B(r_{ m+1})\times (t_0, t)$, $ t_0 < t < 0$, sum over $i=1,2,3$. Then, applying the integration by parts, we have
 \begin{align*}
 & \| \nabla \Omega(t) \phi _m^6 \|_{ q}^{ q}  
 -6q\intl_{t_0}^{t} \intl_{B( r_{ m+1})  \setminus B(r_m)} \frac{W\cdot y }{|y|}  
 | \nabla \Omega(s) |^q \phi _m^{ 6q-1} \eta_m '(| y|) dy ds
 \\
 &\quad \le  \| \nabla \Omega(t_0) \phi_m  ^6\|_{ q}^{q} - \frac{q-3}{2} \intl_{t_0}^{t} \intl_{B(r_{ m+1})}
\frac{ (-s)^{ -  \frac{1}{2} }}{1+ (-s)^{   \frac{1}{2} } } | \nabla \Omega(s) |^q  \phi _m ^{ 6q} dyds 
\\
&\qquad +q\intl_{t_0}^{t} \intl_{B(r_{ m+1})} \frac{1}{1+(-s)^{\frac12}}  \nabla V(s): \partial _i\Omega(s) \otimes \partial _i\Omega(s)  | \nabla \Omega(s) |^{ q-2} \phi_m ^{ 6q} dyds 
\\
&\qquad +q\intl_{t_0}^{t} \intl_{B(r_{ m+1})}   \frac{1}{1+(-s)^{\frac12}} \partial _i \nabla V(s) : \partial _i \Omega(s) \otimes \Omega (s)  | \nabla \Omega(s) |^{ q-2} \phi_m ^{ 6q} dyds. 
 \end{align*}
Since  $ B(r_{ m+1})  \setminus B(r_{ m}) \subset B(\rho _0)  \setminus B(\rho /2) $, the fact  \eqref{7.a2} ensures that 
  $ W\cdot y >0 $ in $ B(r_{ m+1})  \setminus B(r_{ m})\times (t_0 , 0)$. Furthermore, recalling that $ \eta_m ' \le 0$, we see that 
  the sign of the  integral on the left-hand side of the above 
 inequality is non-negative.  Consequently,  taking into account $ q>3$,  it follows that 
 \begin{align}
 & \| \nabla \Omega(t) \phi_m  ^6\|_{ q}^q  
 \cr
 &\quad \le   \| \nabla \Omega(t_0) \phi _m^6 \|_{ q}^q  
+q\intl_{t_0}^{t} \intl_{B(r_{ m+1})}  \frac{1}{1+(-s)^{\frac12}} \nabla V(s): \partial _i\Omega(s) \otimes \partial _i\Omega(s)  | \nabla \Omega(s) |^{ q-2} \phi_m ^{ 6q} dyds 
\cr
&\qquad +q\intl_{t_0}^{t} \intl_{B(r_{ m+1})}  \frac{1}{1+(-s)^{\frac12}}\partial _i \nabla V(s) : \partial _i \Omega(s) \otimes \Omega (s)  | \nabla \Omega(s) |^{ q-2} \phi_m ^{ 6q} dyds
\cr
&\quad =    \| \nabla \Omega(t_0) \phi_m ^6 \|_{ q}^q  +  I+ II.
\label{7.0d}
 \end{align}
 For notational simplicity we define 
\bb\label{zm}
Z_m (s) = \| \nabla^2 V (s) \phi_m ^6 \|_{ q}^{ q},\quad  t_0 \le s < 0.    
\ee
Observing the pointwise estimate $ |\nabla \Omega| \le 2 |D^2 V|$, we immediately see that  
 \begin{align*}
 I+II&\le q2^{ q+1}\intl_{t_0}^{t} \| \nabla V(s)\|_{ L^\infty(B(r_{ m+1}))} Z_m  (s) ds.
 \end{align*}
  
Inserting the estimates of $ I+II$ into  \eqref{7.0d}, we deduce  
\begin{align}
\label{7.0e}
& \| \nabla \Omega  (t) \phi _m ^6\|_{ q}^q
 \le 2^{q}Z_m  (t_0)  + q 2^{ q+1} \intl_{t_0}^{t} \| \nabla V(s)\|_{ L^\infty(B(r_{ m+1}))} Z_m  (s) ds.
\end{align}
Furthermore, in view of \eqref{8.22}, we see that
\begin{align}
\label{7.0dd}
Z_m(t) &\le  c\| \nabla \Omega  (t) \phi _m ^6\|_{ q}^q + c \| \nabla \phi _m \|_{ \infty}^{ \frac{7q-6}{2q}} 
\| V(t) \phi _m\|_{ 2} 
\cr
&\le c\| \nabla \Omega  (t) \phi _m ^6\|_{ q}^q + 
c r^{ - \frac{7q-6}{2q}} \Big(\frac{\rho _0+r}{\rho _0-r}\Big)^{   \frac{(7q-6)m}{2q}}. 
\end{align} 
 Combining \eqref{7.0dd} with \eqref{7.0e}, we find that 
\begin{align}
Z_m (t)   \le cZ_m  (t_0)  + c r^{ - \frac{7q-6}{2q}} \Big(\frac{\rho _0+r}{\rho _0-r}\Big)^{  \frac{(7q-6)m}{2q}} + c \intl_{t_0}^{t} \| \nabla V(s)\|_{ L^\infty(B(r_{ m+1}))} Z_m  (s) ds.
\label{7.0g}
\end{align}

By means of the local version of the logarithmic Sobolev embedding inequality (cf. Lemma\,\ref{lem8.2}) we find for 
every $ s\in (t_0, 0)$
\begin{align}
&\| \nabla V(s) \|_{ L^\infty(B(r_{ m+1}))} 
\cr
&\le c(1+ \| \nabla V(s)\|_{ BMO(B(\rho_{ \ast}))} ) \log (e+ 
\| \nabla^2 V(s)\|_{ L^{ q}(B(r_{ m+1}))} + \| V(s)\|_{ L^2(B(r_{ m+1}))})
\cr
&\le c(1+ \| \nabla V(s)\|_{ BMO(B(\rho_{ \ast}))} ) \log (e+ \| v(s)\|_{ L^2(B(\rho ))}+ 
Z_{ m+1})\cr
&\le c(1+ \| \nabla V(s)\|_{ BMO(B(\rho_{ \ast}))} )\log (e+ 
Z_{ m+1}) ,
\label{7.8}
\end{align}
where the constant $ c>0$ depends only on $ \rho$.   
 
 \hspace{0.5cm}
We continue our discussion by estimating the term on the right-hand side involving the $ BMO$ norm of $ \nabla V$. For this purpose 
let  $ \eta \in  C^{\infty} (\R)$  be a cut off function such that $ 0 \le \eta   \le  1 $  
in $ \R$, $ \eta   \equiv  1$ on $ B(\rho_{ \ast} )$, $ \eta \equiv 0$ in $ (\rho _0, + \infty)$, and $ 0 \le  -\eta  ' \le 2(\rho _0- \rho _{ \ast})^{ -1}= 4 (\rho _0- r)^{ -1}$.   We set $ \psi (y)= \eta (| y|)$. By means of Jensen's inequality we get for every ball $ B(x_0, \sigma )$ with $ x_0\in B(\rho _0/2)$ and   $ 0< \sigma < \rho _0$
\begin{align*}
&\intmw_{B(\rho_{ \ast})\cap B(x_0, \sigma )} | \nabla V(s) - (\nabla V(s))_{B(\rho _0/2)\cap B(x_0, \sigma ) }| dy 
\\
&\quad \le c\intmw_{B(\rho_{ \ast})\cap B(x_0, \sigma )} \intmw_{B(\rho _0/2)\cap B(x_0, \sigma )} | \nabla V(y, s) - \nabla V(y', s)| dy dy' 
\\ 
&\quad \le c\intmw_{ B(x_0, \sigma )} | (\nabla V(s)) \psi^6 - (\nabla V(s) \psi^6 )_{ B(x_0, \sigma ) }| dy. 
\end{align*}
 This together with \eqref{8.15e} with $ r=\rho _0$  yields
 \begin{align}
&\| \nabla V(s)\|_{ BMO(B(\rho_{ \ast}))} 
\cr
& \le  \| \nabla V(s) \psi ^6\|_{ BMO}
\cr
&\le c\Big\{1+ \rho _0^{ 5}\| \nabla \psi \|_{ \infty}^{ 5}\Big\} |\Omega (s) |_{ BMO(B(\rho _0))} +
c\Big\{\rho _0^{ -\frac{5}{2}} + \rho _0^{ \frac{5}{2}}\| \nabla \psi \|^{ 5}_{ \infty}\Big\}
\| V(s) \|_{L^2(B(\rho _0))}
 \cr
&\le c \Big(\frac{\rho _0}{\rho _0- \rho _{ \ast}}\Big)^5|\Omega (s) |_{ BMO(B(\rho _0))} +
c \rho ^{ - \frac{5}{2}}\Big(\frac{\rho _0}{\rho _0- \rho _{ \ast}}\Big)^5
\| V(s) \|_{L^2(B(\rho _0))}
\cr
&\le c \Big(\frac{\rho _0}{\rho _0- \rho _{ \ast}}\Big)^5 \Big\{|\omega (s) |_{ BMO(B(\rho ))} +
c \rho ^{ - \frac{5}{2}}
\| v \|_{L^\infty(-\rho , 0;L^2(B(\rho )))}\Big\}\cr
&\le c( |\omega (s) |_{ BMO(B(\rho ))} +1),
\label{7.9}
 \end{align}
 where we used Corollary \ref{cor10.1} in the last step.
  Combining \eqref{7.8} and \eqref{7.9}, we get 
\begin{align}
&\| \nabla V(s) \|_{ L^\infty(B(r_{ m+1}))} 
\le c\Big\{1+ |\omega (s)|_{ BMO(B(\rho))} \Big\} \log (e+
Z_{ m+1}).
\label{7.9b}
\end{align}   
Combining \eqref{7.0g} and \eqref{7.9b}, we arrive at   
\begin{align}
Z_m (t)   \le d^m  + 
\intl_{t_0}^{t} a(s) Z_m  (s)  \log (e+ Z_{ m+1})ds,
\label{7.0x}
\end{align}
where
\[
  d= c\Big(\frac{\rho _0+r}{\rho _0-r}\Big)^{ \frac{7q-6}{2q}},\quad  a(s) = c(1+ | \omega (s)|_{ BMO(B(\rho ))}), 
\]
while $ c$ denotes a positive constants independently on $ m\in \N$. Setting 
\bb\label{ym}
Y_{ m}(s)= Z_m(s) + e, 
\ee
and eventually replacing $ c$ by a larger constant independent on $ m$, the estimate \eqref{7.0x} turns into   
\begin{align}
e+ Y_m (t)   \le  d^m  + 
\intl_{t_0}^{t} a(s) Y_m  (s)  \log (e+Y_{ m+1}) ds.
\label{7.0gg}
\end{align}
We define 
\bb\label{xm}
X_m(t) = d^m+ \intl_{t_0}^{t} a(s) Y_m  (s)  \log (e+Y_{ m+1}(s))  ds, \quad  t\in [t_0, 0). 
\ee
Then,  \eqref{7.0gg}   implies  $ e+ Y_m \le  X_m$, and thus  the differential inequality 
\begin{equation}
X'_m = a Y_m  \big(\log (Y_{ m+1})+1\big) \le a X_m(t) \log(X_{ m+1}) \quad  \text{ in}\quad  [t_0, 0). 
\label{d_inequ}
\end{equation}
Dividing both sides by $ X_m$, we are led to 
\[
\beta _m ' \le a(t) \beta _{ m+1},\quad  \text{ where}\quad  \beta _m(t) := \log(X_{ m+1}(t)). 
\]
Integrating the both sides over $ (t_0, t)$ with $t_0 \le t < 0$, we find 
\[
\beta _m(t) \le m\log d +  \intl_{t_0}^{t}  a(s) \beta _{ m+1}(s) ds. 
\]
We now verify that the sequence $\{ \beta_m (t)\}$ satisfy the condition \eqref{A.4a} of Lemma \ref{lemA.4} below. 
From the definitions \eqref{xm}, \eqref{ym} and \eqref{zm} one has
\begin{align*}
&\beta_m (t)=\log (X_{m+1} (t))  \\
&\quad\le m\log d +\log\bigg( \int_{t_0} ^{t_1} |a(s)|ds\big) \sup_{t_0<s<t}\log \left\{|Y_m (s)| \log (e+ Y_{m+1}(s))\right\}\bigg) +1\\
&\quad\le  m\log d\\
&\qquad +3 \log\bigg( \int_{t_0} ^{t_1} |a(s)|ds\big) \log \left\{\sup_{t_0<s<t}\| \nabla ^2 V(s)\|_{ L^q(B(\rho_0))}  +e\right\}\bigg)+1\\
&\quad \le M(t)^m
\end{align*}
for all $m\in \Bbb N$, where we set
$$M(t) =\log d + 3 \log\big( \int_{t_0} ^{t_1} |a(s)|ds\big) \log \left\{\sup_{t_0<s<t}\| \nabla ^2 V(s)\|_{ L^q(B(\rho_0))}  +e\right\} +2<+\infty $$
for each $t\in [t_0, t_1)$. Therefore the condition   \eqref{A.4a}  is satisfied.
Applying  Lemma\,\ref{lemA.4},  it follows that 
\begin{align}
\log(e+ Y_0(t))  
 \le  \beta _0(t) \le  \log d  \intl_{t_0}^{t}  a(s)dse^{ \intl_{t_0}^{t}  a(s)ds}. 
\label{7.15}
\end{align}
According to the hypothesis \eqref{bkm} we see that $ \sup_{ t\in (-\rho ,0)}\log(e+ Y_0(t)) <+ \infty $. This yields   
\[
\sup_{ t\in (-\rho ,0)}\| \nabla ^2 v(t)\|_{ L^q(B(r))} \le \sup_{ t\in (-\rho ,0)}\| \nabla ^2 V(t)\|_{ L^q(B(\rho_0))} < +\infty. 
\]
This   completes the proof of Theorem\,\ref{thm1.1}.  \hfill \Beweisende  \\
\ \\
  $$\mbox{\bf Acknowledgements}$$
Chae was partially supported by NRF grants 2016R1A2B3011647, while Wolf has been supported 
supported by the NRF grand 2017R1E1A1A01074536.  
  The authors declare that they have no conflict of interest.    No data sets were generated or analyzed during the current study.
  \appendix

\section{ Gronwall type lemma}
\label{sec:-A}
\setcounter{secnum}{\value{section} \setcounter{equation}{0}
\renewcommand{\theequation}{\mbox{A.\arabic{equation}}}}

\begin{lem}
\label{lemA.2}
Let $ -\infty<a<b<+\infty$. Let $ f\in L^1(a,b)$ and  $ k\in \N$. Then  it holds 
\begin{equation}
 \intl_{a}^{b} \intl_{a}^{t_1} \cdots \intl_{a}^{t_{ k-1}}  \prod_{ j=1}^{ k} f(t_j) d t_k  d t_{ k-1} \ldots d t_1
 = \frac{1}{k!}  \bigg(\intl_{a}^{b} f(t) dt\bigg)^k. 
\label{A.2}
\end{equation}

\end{lem}

{\bf Proof}: We prove the assertion by induction. For $ k=1$ \eqref{A.2} is obvious. Assume the assertion holds for $ k-1$.   
Using the hypothesis of induction, we find 
 \begin{align}
 \intl_{a}^{b} \intl_{a}^{t_1} \cdots \intl_{a}^{t_{ k-1}}  \prod_{ j=1}^{ k} f(t_j) d t_k  d t_{ k-1} \ldots d t_1  &=   \intl_{a}^{b} f(t_1)\bigg\{ \intl_{a}^{t_1} \cdots \intl_{a}^{t_{ k-1}}  \prod_{ j=2}^{ k} f(t_j) d t_k  d t_{ k-1} \ldots dt_2 \bigg\}d t_1
 \cr
 & = \frac{1}{(k-1)!}\intl_{a}^{b} f(s)   \bigg(\intl_{a}^{s} f(t) dt\bigg)^{ k-1} ds.
 \label{A.3}
 \end{align}
Applying integration by parts, we calculate 
\begin{align*}
&\intl_{a}^{b} f(s)   \bigg(\intl_{a}^{s} f(t) dt\bigg)^{ k-1} ds 
 =\intl_{a}^{b} \frac{d}{ds} \intl_{a}^s   f(\tau ) d\tau    \bigg(\intl_{a}^{s} f(t) dt\bigg)^{ k-1} ds
\\
&\qquad =  \bigg(\intl_{a}^{b} f(t) dt\bigg)^k - (k-1) \intl_{a}^{b} f(s) \intl_{a}^s   f(\tau ) d\tau     \bigg(\intl_{a}^{s} f(t) dt\bigg)^{ k-2} ds
\\
&\qquad =  \bigg(\intl_{a}^{b} f(t) dt\bigg)^k - (k-1) \intl_{a}^{b} f(s)      \bigg(\intl_{a}^{s} f(t) dt\bigg)^{ k-1} ds. 
\end{align*}
This yields
\begin{equation}
\intl_{a}^{b} f(s)   \bigg(\intl_{a}^{s} f(t) dt\bigg)^{ k-1} ds = \frac{1}{k} \bigg(\intl_{a}^{b} f(t) dt\bigg)^k. 
\label{A.4}
\end{equation}
Replacing the integral on the right-hand side of \eqref{A.3} by the  right-hand side of \eqref{A.4}, we obtain \eqref{A.2} for $ k$. 
Hence by induction \eqref{A.2} holds for all $ k\in \N$.  \hfill \Beweisende

 \begin{lem}[Iteration lemma]
\label{lemA.4}
Let $ \beta  _m: [t_0, t_1] \rightarrow \R$, $ m\in \N\cup \{ 0\}$ be a sequences of bounded functions. Suppose there exists $0<K=K(t) <+\infty $ for each $t\in [t_0, t_1)$ such that
\bb\label{A.4a}
  |\beta_m(t)| < K(t)^m \qquad \forall t\in [t_0, t_1), \forall m\in \Bbb N.
\ee

Furthermore let 
$ a\in L^1(t_0, t_1)$ with $a(t)\geq 0$ for almost every $t\in [t_0, t_1]$. We assume that the following recursive  integral 
inequality holds true for a constant $ C>0$
\begin{align}
\beta _{ m}(t) \le  Cm +   \intl_{t_0}^{t} a(s) \beta _{ m+1}(s)  ds, \quad  m\in \N\cup \{ 0\}.
\label{A.5a}
\end{align}
Then the following inequality holds true for all $ t\in [t_0,t_1]$
\begin{equation}
\beta _0(t) \le  C  \intl_{t_0}^{t} a(s)  ds \, e^{\intl_{t_0}^{t} a(s )   ds }. 
\label{A.5b}
\end{equation}
\end{lem}

{\bf Proof}:  Iterating \eqref{A.5a} $ m$-times, and applying Lemma \ref{lemA.2}, we see that for each $t\in [t_0, t_1)$ it follows
\begin{align}
\beta _0(t) &\le    C\intl_{t_0}^{t} a(s_1)   ds_1 + 2C\intl_{t_0}^{t} \intl^{s_1}_{t_0} a(s_1) a(s_2)   ds_2 ds_1 
\cr
&\qquad + \ldots + Cm\intl_{t_0}^{t} \intl^{s_1}_{t_0} \ldots \intl^{s_{ m-1}}_{t_0} a(s_1) a(s_2) \ldots a(s_m)ds_{ m} \ldots ds_2 ds_1 
\cr
&\qquad + \intl_{t_0}^{t} \intl^{s_1}_{t_0} \ldots \intl^{s_{ m}}_{t_0} a(s_1) a(s_2) \ldots a(s_{ m+1}) \beta _{ m+1}(s_{ m+1})  
ds_{ m+1} \ldots ds_2 ds_1 
\cr
&= \sum_{k=1}^{m}\frac{C}{(k-1) !} \bigg(\intl_{t_0}^{t}  a(s)   ds \bigg)^{ k}+ J_m(t),
\label{A.14}
\end{align} 
where
\begin{align*}
|J_m (t)|&\le  \sup_{0<s< t} |\beta_{m+1}(s)| \intl_{t_0}^{t} \intl^{s_1}_{t_0} \ldots \intl^{s_{ m}}_{t_0} a(s_1)  a(s_2) \ldots a(s_{ m+1})  
ds_{ m+1} \ldots ds_2 ds_1\\
&\le \frac{ 1}{(m+1)!}  \bigg(\sup_{0<s< t} K(s)  \intl_{t_0}^{t}  a(s )   ds\bigg)^{ m+1} \to 0\quad \text{as $\quad m\to +\infty$}
\end{align*}
for each $t\in [t_0, t_1)$.
Therefore, 
\begin{align*}
\beta_0(t)\leq \sum_{k=1} ^\infty  \frac{C}{(k-1) !} \bigg(\intl_{t_0}^{t}  a(s)   ds \bigg)^{ k}= C  \intl_{t_0}^{t} a(s)  ds\,  e^{\intl_{t_0}^{t} a(s )   ds}. 
\end{align*}
 \hfill \Beweisende

\section{Gagliardo-Nirenberg's inequality with cut-off}
\label{sec:-B}
\setcounter{secnum}{\value{section} \setcounter{equation}{0}
\renewcommand{\theequation}{\mbox{B.\arabic{equation}}}}

\begin{lem}
\label{lemB.1}
Let $ \psi \in C^{\infty}_{\rm c}(\R^{n})$ with $ 0 \le \psi \le 1$, and $k \ge 1 $.  For all $ u\in W^{1,\, q}(\R^{n})\cap L^2(\R^{n})$, $ 2 < q < +\infty$, and 
for $ m \ge k\frac{n(q-2)+ 2q}{2q}$  it holds 
\begin{align}
&\| u \psi ^{ m-k}\|_{ q}
\cr
 &\le  c \| u \psi ^{ m - k\frac{n(q-2)+ 2q}{2q}}\|_{ 2}^{ \frac{2q}{n(q-2)+ 2q}}\| \nabla  u \psi ^m\|_{ q}^{ \frac{n(q-2)}{n(q-2)+ 2q}} + c\| \nabla \psi \|_{ \infty}^{\frac{n(q-2)}{2q}}\| u \psi ^{ m - k\frac{n(q-2)+ 2q}{2q}}\|_{ 2} . 
\label{B.1}
\end{align}
\end{lem}

{\bf Proof}:  Let $ \psi \in C^{\infty}_{\rm c}( \R^{n})$ with $ 0 \le \psi \le 1$, and $ k \ge 1$. Let  $  m \ge k\frac{n(q-2)+ 2q}{2q}$. By means of H\"older's inequality we estimate 
\begin{align*}
\| u \psi^{ m-k} \|_{ q} &\le \| u \psi ^{m - k- k\frac{n(q-2)}{2q} }\|_{ 2}^{ \frac{2}{n(q-2)+2}} 
\| u \psi^{m- k + \frac{k}{q}} \|_{ \frac{nq}{n-1}}^{ \frac{n(q-2)}{n(q-2)+2}}
\\
&\le c\| u \psi ^{m - k- k\frac{n(q-2)}{2q} }\|_{ 2}^{ \frac{2}{n(q-2)+2}} 
\|w\|_{ \frac{n}{n-1}}^{ \frac{n(q-2)}{qn(q-2)+2q}}, 
\end{align*}
where $ w = | u|^q \psi^{qm- kq + k}$.
By  Sobolev's inequality  along with H\"older's inequality we infer
\begin{align*}
\|w\|_{ \frac{n}{n-1}} &\le \| \nabla w\|_{ 1}
\\
&\le c\intl_{ \R^{n}}  | u|^{ q-1}| \nabla u| \psi^{qm- kq + k} dx + c \| \nabla \psi \|_{ \infty}\| u \psi ^{ m-k}\|^q_{ q}
\\
&\le c\| u \psi ^{ m-k}\|_{ q}^{ q-1} \| \nabla u \psi ^m\|_{ q}  + c \| \nabla \psi \|_{ \infty}\| u \psi ^{ m-k}\|^q_{ q}.
\end{align*}
Combining the last two estimates and applying Young's inequality, we obtain the assertion \eqref{B.1}.  
 \hfill \Beweisende 
 
 \begin{lem}
\label{lemB.2}
Let $ u \in W^{1,\, q}(B(\rho ))$, $ n< q < +\infty$ such that  $ \nabla \cdot u =0$ almost everywhere in $ B(\rho )$. 
Then for every   $ \psi \in C^{\infty}_c(B(\rho ))$ with $ 0 \le \psi \le 1$ and $ m \ge \frac{q}{2}$  it holds  
\begin{equation}
\| u \psi ^{ m}\|_{ \infty} \le c\| u \psi ^{ m- \frac{q}{2}}\|_{ 2}^{ \frac{2(q-n)}{2q-2n+nq } }   \| \nabla  u \psi ^{ m}\|^
{ \frac{nq}{2q-2n+nq}}_{ q} + 
c\| \nabla \psi \|_{ \infty}^{ \frac{n}{2}} \| u \psi ^{ m-\frac{q}{2}}\|_{ 2}. 
\label{B.2}
\end{equation}
\end{lem} 

{\bf Proof}: By virtue of Gagliardo-Nirenberg's inequality we estimate  
\begin{align*}
\| u \psi ^{ m}\|_{ \infty} &\le c \| u \psi ^{ m}\|_{ q}^{ 1- \frac{n}{q}} \|\nabla ( u \psi ^{ m})\|_{ q}^{ \frac{n}{q}}
\\
&\le c\| u \psi ^{ m-1} \|_{ q}^{ 1- \frac{n}{q}}   \| \nabla u \psi ^{ m}\|^{   \frac{n}{q}}_{ q} + 
c\| \nabla \psi \|_{ \infty}^{ \frac{n}{q}} \| u \psi ^{ m-1}\|_{ q}.
\end{align*}
From $L^q$-interpolation  we find
\[
\| u \psi ^{ m-1}\|_{ q} \le \| u \psi ^{ m- \frac{q}{2}} \|_{ 2}^{ \frac{2}{q }} \| u \psi ^{ m}\|_{ \infty}^{ 1- \frac{2}{q}}. 
\]
Combining the  two estimates above, and applying Young's inequality, we obtain the assertion of the lemma.  \hfill \Beweisende 

\vspace{0.2cm}
Although the following result of  embedding  of $ BMO(B(r))$  into $ L^p(B(r))$  might be well known, for readers convenience we give a short proof.  

\begin{lem}
\label{lemB.3}
Let $ u\in BMO(B(r))$. Then $ u\in \cap_{ 1 \le  q< \infty} L^q(B(r))$, and it holds 
\begin{equation}
\| u\|_{ L^q(B(r))} \le c r^{ \frac{n}{q}} | u|_{BMO(B(r))}+ c r ^{ \frac{n}{q}-n}\| u\|_{ L^1(B(r))} =
c r^{ \frac{n}{q}} \| u\|_{ BMO(B(r))}. 
\label{B.3}
\end{equation}
\end{lem}

{\bf Proof}: 
Let $ u\in BMO(B(1))$. By $ U\in BMO$ we denote the extension defined in Section\,2.  By John-Nirenberg's inequality\cite{joh} it holds 
\[
m(\{x\in Q(2) \,|\, | U - U_{ Q(2)} | > \lambda \}) \le c_1  e^{ - \frac{c_2 \lambda}{| U| _{ BMO(Q(2))}}}\quad \forall\,\lambda\in (0, +\infty)
\]
with constants $ c_1, c_2$, depending on $ n$ and $ q$ only, where $m(\cdot)$ denotes the Lebesgue measure in $\Bbb R^n$, and 
$Q(r)$ denotes the cube with the center  at origin and the side length $2r$.
Multiplying both sides by $ q\lambda^{ q-1}$, integrating the result over $ (0, +\infty)$, using a suitable change of coordinates, and 
employing \eqref{8.2a},   we arrive at 
\begin{align*}
\| U- U _{ Q(2)} \|_{ L^q(Q(2))}^q &= 
q  \intl_{0}^{\infty}   \lambda^{ q-1}m(\{x\in Q(2) \,|\, | U- U _{ Q(2)} | > \lambda \} )
\\
&\le q | U| _{ BMO(Q(2))}^q  \intl_{0}^{\infty}   \lambda^{ q-1}c_1  e^{ - c_2 \lambda }
\\
&\le c \| u\|^q_{ BMO(B(1))}. 
\end{align*}
Accordingly, 
\[
\| u\|_{ L^q(B(1))} \le \|U- U_{ Q(2)} \|_{ L^q(Q(2))}+ 
c| U_{ Q(2)}| \le c \| u\|_{ BMO(B(1))}. 
\]
Hence, \eqref{B.3} follows from the above estimate by using  a standard scaling  argument.  \hfill \Beweisende

\bibliographystyle{siam} 

\begin{thebibliography}{99}
 \bibitem{bar} C. Bardos and E. S.  Titi, {\it Euler equations for an ideal incompressible fluid,}  Russian Math. Surveys,  {\bf 62}, no. 3,  (2007),  pp. 409-451 \bibitem{bea}J. T. Beale, T. Kato and A. Majda,  {\it Remarks on the
breakdown of smooth solutions for the 3-D Euler equations}, Comm.
Math. Phys., {\bf 94},  (1984), pp. 61-66.
\bibitem{cha}D. Chae, {\it  Incompressible Euler equations: the blow-up problem and related results,} Handbook of differential equations: evolutionary equations. Vol. IV,   Elsevier/North-Holland, Amsterdam, (2008), pp. 1-55.
\bibitem{cha1}D. Chae and J. Wolf, {\it On the local Type I conditions for the 3D Euler equations,}  
Arch. Rat. Mech. Anal., (2) (2018), pp 641-663.
\bibitem{con}P. Constantin,  {\it On the Euler equations of incompressible fluids,} Bull. Amer. Math. Soc., {\bf
44}, no. 4, (2007), pp. 603-621.
 \bibitem{con1}P. Constantin, C. Fefferman and A. Majda,
 {\it Geometric constraints on potential singularity formulation in the
 3-D Euler equations}, Comm. P.D.E., {\bf 21}, (3-4), (1996), pp.
 559-571. 
   \bibitem{den}J. Deng, T. Y. Hou and X. Yu, {\it Improved geometric
conditions for non-blow up of the 3D incompressible Euler equations,}
Comm. P.D.E., {\bf 31}, no. 1-3, (2006), pp. 293-306.
 \bibitem{eul}L. Euler, {\it Principes g\'{e}n\'{e}raux du mouvement des fluides,}
M\'{e}moires de l'acad\'{e}mie des sciences de Berlin, {\bf 11},
(1755), pp. 274-315. 
\bibitem{joh}F. John and L. Nirenberg, {\it On functions of bounded mean oscillation}, Comm. Pure Appl. Math. {\bf 343}, (1983), pp. 415-426.
\bibitem{kat} T. Kato and  G. Ponce, {\it Commutator estimates and the Euler and Navier-Stokes equations,}  Comm. Pure Appl. Math., {\bf 41} no. 7,(1988),  pp. 891-907.
  \bibitem{ker} R. M. Kerr, {\it  Evidence for a singularity of the three-dimensional, incompressible Euler equations,}
   Phys. Fluids,  A 5 (1993), no. 7, pp. 1725-1746.
    \bibitem{koz1}
  H. Kozono and Y. Taniuchi, {\it  Bilinear estimates in BMO and the Navier-Stokes equations,}
  Math. Z. {\bf 235},  (2000), no. 1, pp. 173-194. 
  \bibitem{koz}H. Kozono and Y. Taniuchi, {\it Limiting case of the
Sobolev inequality in BMO, with applications to the Euler
equations,} Comm. Math. Phys., {\bf 214}, (2000), pp. 191-200.   
   
 \bibitem{hou}G. Luo and T. Hou, {\it  Toward the finite-time blowup of the 3D axisymmetric Euler equations: a numerical investigation},
 Multiscale Model. Simul. {\bf 12} (2014), no. 4, pp. 1722-1776. 

\bibitem{maj}A. Majda and A. Bertozzi, {\it Vorticity and
Incompressible Flow,} Cambridge Univ. Press. (2002).
\bibitem{ste}E. M. Stein, {\it Harmonic Analysis:  Real-variable Methods, Orthogonality, and Oscillatory Integrals, } Princeton University Press, (1993).
\bibitem{tao}T. Tao, {\it  Finite Time Blowup for Lagrangian Modifications of the Three-Dimensional Euler Equation}, Ann. PDE, {\bf 2}, no. 9, (2016), pp. 1-79.
\end{thebibliography}

\end{document}